\documentclass[a4paper,10pt]{article}

\usepackage{amsfonts}
\usepackage{amssymb}
\usepackage{amsmath}
\usepackage{amsthm}
\usepackage{enumerate}
\usepackage{epsfig}
\usepackage{mathrsfs}
\usepackage{graphicx}
\usepackage{dsfont}
\usepackage{epsfig,color}
\usepackage{dsfont}
\usepackage[utf8]{inputenc}
\usepackage{bm}
\usepackage{titlesec}
\titleformat{\section}[block]{\Large\scshape\filcenter}{\thesection}{1em}{}

\newcommand{\R}{\mathbb{R}}

\newcommand{\W}{\mathcal{W}}

\newcommand{\Q}{{[0,1)^d}}

\newcommand{\Td}{{\mathbb{T}^d}}

\newcommand{\meni}{\leqslant}

\newcommand{\PT}{\mathcal{P}(\Td)}

\newtheorem{definicion}{Definition}[section]
\newtheorem{definition}[definicion]{Definition}
\newtheorem{proposition}[definicion]{Proposition}
\newtheorem{theorem}[definicion]{Theorem}

\newtheorem{lemma}[definicion]{Lemma}
\newtheorem{remark}[definicion]{Remark}

\numberwithin{equation}{section}

\def\fin { \vskip 0pt \hfill \hbox{\vrule height 5pt width 5pt depth 0pt} \vskip 12pt}

\begin{document}

\title{{\large\textbf{MINIMIZING MOVEMENT FOR A FRACTIONAL POROUS  \\ MEDIUM EQUATION IN A PERIODIC SETTING}}}
\date{}

\author{\textsc{\small Lucas C. F. Ferreira,\hspace{.2cm} Matheus C. Santos \hspace{0.2cm}and\hspace{0.2cm} Julio C. Valencia-Guevara} }

\maketitle

\footnotetext{Departamento de Matem\'atica - IMECC, Universidade Estadual de Campinas, 13083-859, Campinas-SP, Brazil.}
\footnotetext{Email addresses: lcff@ime.unicamp.br (LCFF), msantos.ufrgs@gmail.com (MCS), ra099814@ime.unicamp.br (JCV-G)}

\footnotetext{\today{}}

\begin{abstract}
	We consider a fractional porous medium equation that extends the classical
	porous medium and fractional heat equations. The flow is studied in the space
	of periodic probability measures endowed with a non-local transportation
	distance constructed in the spirit of the Benamou-Brenier formula. For initial
	periodic probability measures, we show the existence of absolutely continuous
	curves that are generalized minimizing movements associated to R\'{e}nyi
	entropy. For that, we need to obtain entropy and distance properties and to
	develop a subdifferential calculus in our setting. \newline
	
%
	
\end{abstract}






\section{\sc Introduction}

\quad\ We are concerned with the fractional porous medium equation (FPME)
\begin{equation}
	\partial_{t}\rho+(-\Delta)^{\sigma}\rho^{m}=0\text{ and }\rho(0,x)=\rho
	_{0}(x)\label{FPMeq}%
\end{equation}
where the spatial-dimension $d\geq1,$ $0<\sigma<1$ and $m_{\ast}%
:=\frac{(d-2\sigma)_{+}}{d}<m\leqslant2.$ Our aim is to study the global flow
$\rho=\rho_{t}(x)$ with spatial-periodic conditions, i.e. $(t,x)\in
\lbrack0,+\infty)\times{\mathbb{T}^{d}}$ (${\mathbb{T}^{d}}$ is the
$d$-dimensional torus). We use an optimal transport approach on the space
$\mathcal{P}({\mathbb{T}^{d}})$ equipped with a pseudo-metric $\mathcal{W}$ to
find a generalized minimizing movement curve associated to the R\'{e}nyi entropy, which gives us a weak notion of gradient flow for \eqref{FPMeq}.

The equation \eqref{FPMeq} appears in the study of hydrodynamic limits of
particle systems \cite{Jara}, the boundary heat control problems as studied by
Athanasopoulos and Caffarelli \cite{athanas-caffar}, and the so-called
anomalous diffusion process that appears when jump processes are introduced in
the modeling (see \cite{AbeThurner,VIKH,WZ}). In
\cite{de-Pablo-Vazquez,de-Pablo-Vazquez2}, the authors studied the model
\eqref{FPMeq} and developed a theory for the problem in $\mathbb{R}^{d}$ and
with general initial data $\rho_{0}\in L^{1}(\mathbb{R}^{d})$ that includes
also sign changing solutions. For the case $\sigma=1/2$ they constructed a
weak solution by rewriting the non-local problem in a local way by mean of a
well known harmonic extension. The general case was treated in
\cite{de-Pablo-Vazquez2} by approximating the initial data by bounded
functions and approximating $\mathbb{R}^{d}$ by bounded domains. Existence and
uniqueness results are true for any positive value $m$ of the nonlinearity if the initial data is, for example, also bounded, but the full theory only
works for values of $m$ larger than the critical value $m_{\ast}$. In
this case, if we take a nonnegative integrable initial data $\rho_{0}$, then
conservation of mass, positivity and regularity hold. Let us observe that, for
the sake of coherence, we are assuming $m_{\ast}<m\leqslant2$ in order to
assure mass and sign conservation for \eqref{FPMeq}, however most of arguments
and results in this paper hold true for the full range $m\in(0,2].$

The optimal transport approach was used by \cite{Otto} in the analysis of
solutions to the classical porous medium equation
\begin{equation}
	\partial_{t}\rho=\Delta\rho^{m}=\nabla\cdot(\rho\nabla U^{\prime}%
	(\rho)),\;\;\;\;\;U(\rho)=\frac{1}{m-1}\rho^{m}.\label{classicalPME}%
\end{equation}
It was showed that solutions to this equation can be interpreted as a gradient
flow $t\mapsto\rho_{t}(x)dx$ of probabilities on $\mathcal{P}(\mathbb{R}^{d})$
associated to the R\'{e}nyi entropy
\[
\mathcal{U}_{m}(\rho):=\frac{1}{m-1}\int_{\mathbb{R}^{d}}\rho(x)^{m}\;dx
\]
with respect to (w.r.t. for short) the $L^{2}$-Wasserstein distance $W_{2}$,
which is defined by
\[
W_{2}(\mu_{0},\mu_{1}):=\inf_{\pi\in\Gamma(\mu_{0},\mu_{1})}%
\left(\int_{{\mathbb{R}^{d}}\times{\mathbb{R}^{d}}}|x-y|^{2}\;d\pi(x,y)\right)^{1/2},\text{ for
	all }\mu_{0},\mu_{1}\in\mathcal{P}(\mathbb{R}^{d}),
\]
where $\Gamma(\mu_{0},\mu_{1})$ is the set of probability measures $\pi$ on
${\mathbb{R}^{d}}\times{\mathbb{R}^{d}}$ with marginals $\mu_{0}$ and $\mu
_{1}$. This method has allowed one
to obtain uniqueness, well-posedness and asymptotic behavior of solutions
(such as contraction property and rates for the entropy decay) to the porous
medium and to a number of other equations (see e.g. \cite{Agueh}%
,\cite{villani},\cite{ambrosiogiglisavare},\cite{JKO}). See also
\cite{CarrilloSlepcev}, \cite{GT} and \cite{VF} for periodic solutions of PDEs
without fractional Laplacian. Besides that, the use of the entropy and the
optimal transport distance also have provided a powerful method to prove
existence of solutions by the so-called Minimizing Movement Scheme (in the
general metric framework as considered in \cite{ambrosiogiglisavare}) or
JKO-scheme (in the Wasserstein space $\mathcal{P}_{2}(\mathbb{R}^{d})$ used in
\cite{JKO,ambrosiogiglisavare,CFFLS,FF,BL} and many others). In the metric
approach, that is obtained in the following way (see
\cite{ambrosiogiglisavare}): if $(X,d)$ is a complete metric space and we
consider a lower semicontinuous functional $\mathcal{E}:X\rightarrow
(-\infty,\infty]$, then for an initial data $u_{0}\in X$ and time step
$\tau>0$, we can consider the implicit discrete scheme defined by
\begin{equation}
	\left\{
	\begin{array}
		[c]{l}%
		u^{0}:=u_{0}\\
		u^{n}\in\text{arg min}\left(  \mathcal{E}(u)+\displaystyle\frac{d^{2}%
			(u,u^{n-1})}{2\tau}\right)  ,\;\;n\geqslant1
	\end{array}
	\right.  \label{eq:discr-schem}%
\end{equation}
that is, $u^{n}$ is any minimizer for the functional $u\mapsto\mathcal{E}%
(u)+\frac{d^{2}(u,u^{n-1})}{2\tau}$. Defining the piecewise constant curve
$u_{\tau}(t)=u^{n}$ if $t\in\lbrack n\tau,(n+1)\tau)$, then a limit curve
$u:[0,\infty)\rightarrow X$ of $u_{\tau}$ as $\tau\rightarrow0$ (up to a
subsequence) is called a \textit{generalized minimizing movement }curve (GMM
curve) starting from $u_{0}$ and it is a weak notion of gradient flow for
$\mathcal{E}$.

A slightly different notion of gradient flow in metric spaces involves the
concept of curves of maximal slope w.r.t. upper gradient (see
\cite{ambrosiogiglisavare}): we say that a function $g$ defined on the metric
space is an upper gradient for the entropy $\mathcal{E}$ if, for all
absolutely continuous curve $v:[0,\infty)\mapsto X$, one has $|(\mathcal{E}%
\circ v)^{\prime}(t)|\leqslant g(v(t))|v^{\prime}|(t)$, where $|v^{\prime
}|(t)$ denotes the metric derivative. And we say that an absolutely continuous
curve $u:[0,\infty)\mapsto X$ is a curve of maximal slope for $\mathcal{E}$
w.r.t the upper gradient $g$ if
\begin{equation}
	\frac{d}{dt}\mathcal{E}(u(t))\leq-\frac{g(u(t))^{2}}{2}-\frac{|u^{\prime}%
		|^{2}(t)}{2}\;,\;\;\mbox{ for almost }t\in\lbrack0,\infty).
	\label{eq:maximal-slope-curves}%
\end{equation}
Let us comment that if we consider the Wasserstein space $(\mathcal{P}%
_{2}(\mathbb{R}^{d}),W_{2})$, then the inequality
\eqref{eq:maximal-slope-curves} becomes essentially equivalent to the
continuity equation of the curve $u(t)$. In the most of cases this continuity
equation becomes the partial differential equation that gives birth to the
functional considered. Another way to explain the relation between the
equation \eqref{classicalPME} and the Wasserstein distance (specially in the
case where the velocity field is a gradient vector field) can be seen by the
dynamic characterization of $W_{2}$, which is known as Benamou-Brenier formula
(see \cite{BenamouBrenier}):
\[
W_{2}(\bar{\mu}_{0},\bar{\mu}_{1})^{2}=\inf_{(\mu,v)\in
	\mathcal{CE}_{1}(\bar{\mu}_{0},\bar{\mu}_{1})}\left\{  \int_{0}^{1}%
\int_{{\mathbb{R}^{d}}}|\nabla v_{t}(x)|^{2}\;d\mu_{t}(x)dt\right\}
\]
where $\mathcal{CE}_{1}(\bar{\mu}_{0},\bar{\mu}_{1})$ is the set of all
sufficiently smooth pairs $(\mu_{t},v_{t})_{t\in[0,1]}$ such that
\[
\left\{
\begin{array}
[c]{l}%
\partial_{t}\mu_{t}+\nabla\cdot(\mu_{t}\nabla v_{t}%
)=0\;,\;\;\mbox{ in }(0,1)\times{\mathbb{R}^{d}}\\
\mu_{0}=\bar{\mu}_{0}\;,\;\;\;\;\;\;\;\mu_{1}=\bar{\mu}_{1}%
\end{array}
\right.
\]

As in \cite{DNS,erbar}, we use these ideas to define a metric, in a
Benamou-Brenier style, on a subset of $\mathcal{P}({\mathbb{T}^{d}})$ that
incorporates the non-local character of the problem \eqref{FPMeq} and we use
it to construct a gradient flow associated to the R\'{e}nyi entropy. For that
matter, we need to obtain properties for the non-local metric and entropy, and
develop a subdifferential calculus in our setting.

In \cite{erbar}, the author showed that solutions to the fractional linear
heat equation (i.e. $m=1$)
\begin{equation}
	\partial_{t}\rho+(-\Delta)^{\sigma}\rho=0\text{ in }(0,\infty)\times
	{\mathbb{R}^{d}}\text{ and }\rho(0,x)=\rho_{0}\text{ } \label{frac-heat-1}%
\end{equation}
can be seen as a gradient flow in a non-local metric (built from a L\'{e}vy
jump kernel) by using, among others, that \eqref{frac-heat-1} generates a
well-behaved semigroup. Due to the presence of $(-\Delta)^{\sigma}$ in
\eqref{FPMeq}, we are led to consider some ideas from \cite{erbar}. However,
we adopt a periodic setting that allows us to obtain the necessary compactness
and convergence of the discrete scheme \eqref{eq:discr-schem}. So, we are
moving in the opposite direction from that of \cite{erbar} in the sense that
we are constructing a flow for \eqref{FPMeq} via the scheme
\eqref{eq:discr-schem}, and not using previously known information about
existence of solutions.

Let us just mention that we consider the $d$-dimensional torus ${\mathbb{T}%
	^{d}}$ as the cube ${[0,1)^{d}}$ equipped with the metric $d_{{\mathbb{T}^{d}%
	}}$ defined by
	\begin{equation}
		d_{{\mathbb{T}^{d}}}(x,y)=\min_{k\in\mathbb{Z}^{d}}|x-y+k|.
		\label{metric-torus-1}%
	\end{equation}

	From \cite{RS} we know that, if $f:{\mathbb{T}^{d}}\rightarrow\mathbb{R}$ is
	sufficiently smooth, the following pointwise formula holds for the fractional
	Laplacian on the torus:
	\[
	(-\Delta)^{\sigma}f(x)=p.v.\int_{{[0,1)^{d}}}\left(  f(x)-f(y)\right)
	K^{\sigma}(x-y)\;dy
	\]
	where the kernel $K^{\sigma}$ is $\mathbb{Z}^{d}$-periodic and defined by
	\[
	K^{\sigma}(x):=C_{d,\sigma}\sum_{k\in\mathbb{Z}^{d}}\frac{1}{|x+k|^{d+2\sigma
		}}\;,\;\;\mbox{ with }\;C_{d,\sigma}:=\frac{4^{\sigma}\Gamma\left(  \frac
		{d}{2}+\sigma\right)  }{\pi^{d/2}\left\vert \Gamma(-\sigma)\right\vert }.
	\]
	Motivated by that, let us call $\rho:(0,\infty)\times\mathbb{T}^{d}%
	\rightarrow\mathbb{R}$ (in some suitable class) a weak solution to
	\eqref{FPMeq} if, for all $\varphi\in C_{c}^{\infty}((0,T);C^{\infty
	}({\mathbb{T}^{d}}))$, one has
	\begin{equation}
		\int_{0}^{T}\int_{{[0,1)^{d}}}\partial_{t}\varphi_{t}(x)\rho_{t}%
		\;dxdt-\frac{1}{2}\int_{0}^{T}\int_{G({[0,1)^{d}})}(\varphi_{t}(x)-\varphi
		_{t}(y))(\rho_{t}(x)^{m}-\rho_{t}(y)^{m})K^{\sigma}(x-y)\;dx\;dydt=0
		\label{eq.pfpmeW}%
	\end{equation}
	where
	\begin{equation}\label{G01}
	G({[0,1)^{d}}):=\{(x,y)\in{[0,1)^{d}}\times{[0,1)^{d}}\;|\;x\neq y\}.
	\end{equation}
	Throughout the paper, we denote by $\overline{\nabla}f$ the discrete gradient
	\begin{equation}
		\overline{\nabla}f(x,y):=f(y)-f(x),\text{ for all }f:\mathbb{R}^{d}%
		\rightarrow\mathbb{R}. \label{discret-grad-1}%
	\end{equation}
	Thus, if $m\neq1$, we can rewrite \eqref{eq.pfpmeW} as
	\begin{align}
		\int_{0}^{T}\int_{{[0,1)^{d}}}\partial_{t}\varphi_{t}(x)\rho_{t}\;dxdt  &
		=\frac{1}{2}\int_{0}^{T}\int_{G({[0,1)^{d}})}\overline{\nabla}\varphi
		_{t}\overline{\nabla}\rho_{t}^{m}K^{\sigma}(x-y)\;dxdydt\nonumber\\
		&  =\frac{1}{2}\int_{0}^{T}\int_{G({[0,1)^{d}})}\overline{\nabla}\varphi
		_{t}\overline{\nabla}\rho_{t}^{m-1}\frac{\overline{\nabla}\rho_{t}^{m}%
		}{\overline{\nabla}\rho_{t}^{m-1}}K^{\sigma}(x-y)\;dxdydt\nonumber\\
		&  =\frac{1}{2}\int_{0}^{T}\int_{G({[0,1)^{d}})}\overline{\nabla}\varphi
		_{t}\overline{\nabla}U_{m}^{\prime}(\rho_{t})\widehat{\rho}_{m}K^{\sigma
		}(x-y)\;dxdydt \label{eq.FMPEweak}%
	\end{align}
	where
	\[
	U_{m}(s):=\frac{s^{m}}{m-1}\;\mbox{ and }\;\;\widehat{\rho}_{m}(x,y):=\frac
	{m-1}{m}\frac{\overline{\nabla}\rho^{m}(x,y)}{\overline{\nabla}\rho
		^{m-1}(x,y)}\text{ }.
	\]
	If $m=1$ we have the same representation if we define
	\[
	U_{1}(s):=s\log{s}\;\mbox{ and }\;\;\widehat{\rho}_{1}(x,y):=\frac
	{\overline{\nabla}\rho(x,y)}{\overline{\nabla}\log{\rho}(x,y)}\;.
	\]
	Therefore, we can write \eqref{eq.FMPEweak} formally as
	\[
	\partial_{t}\rho-\overline{\nabla}\cdot(\widehat{\rho_{m}}\overline{\nabla
	}U_{m}^{\prime}(\rho))=0
	\]
	which resembles the classical porous medium equation given in \eqref{classicalPME}.
	
	This paper is organized as follows. In Section 2 we give some basic
	definitions and notations. Section 3 is devoted to the action functional
	$\mathcal{A}(\mu,\bm{\nu})$ and its properties. Some results about a
	periodic non-local continuity equation are presented in Section 4. In Section
	5 we define the periodic non-local Wasserstein distance $\mathcal{W}$ based on
	the Benamou-Brenier formula and show some of its properties. In Section 6 a
	subdifferential calculus is constructed in the periodic Wasserstein space
	$\mathcal{P}(\mathbb{T}^{d})$ endowed with the non-local metric $\mathcal{W}$.
	In Section 7 we show the existence of GMM curves associated to \eqref{FPMeq}.
	

	\section{Notation and preliminaries}
	
	\quad In this section, we make some definitions and remarks about the notation that
	will be used in the sequel.
	
	Following \cite{GT}, we define the equivalence relation in $\mathcal{P}%
	(\mathbb{R}^{d})$ given as
	\begin{equation*}
		\mu_1,\mu_2\in\mathcal{P}(\mathbb{R}^{d}). \;\;\;\mu_1\sim\mu_2\equiv\int
		_{{\mathbb{R}^{d}}}\zeta(x)\; d\mu_1(x) = \int_{{\mathbb{R}^{d}}}\zeta(x)\;
		d\mu_2(x)\;, \forall\;\zeta\in C({\mathbb{T}^{d}}), 
	\end{equation*}
	where $C(\mathbb{T}^{d})$ stands for the space of continuous functions in
	$\mathbb{R}^{d}$ which are $\mathbb{Z}^{d}$-periodic. Then, it is well known
	that $\mathcal{P}({\mathbb{T}^{d}}) = \;\; \mathcal{P}(\mathbb{R}^{d})/\sim$.
	
	Considering the set
	\[
	G({\mathbb{R}^{d}})=\{(x,y)\in{\mathbb{R}^{d}}\times{\mathbb{R}^{d}}\;|\;x\neq
	y\}
	\]
	and the space of functions
	\[
	C_{c}(G(\mathbb{T}^{d}))=\{\varphi\in C(G({\mathbb{R}^{d}}))\;|\;\min
	\{|x-y|\;|\;(x,y)\in\text{supp }(\varphi)\}>0\;,\;\varphi
	\;\mbox{is}\;(\mathbb{Z}^{d}\times\mathbb{Z}^{d})\text{-periodic}\}\;,
	\]
	we can analogously define the equivalence relation
	\[
	\bm{\nu}_1,\bm{\nu}_2\in\mathcal{M}_{0}(\mathbb{R}^{d})\;\;\;\bm{\nu
	}_1\bm{\sim}\bm{\nu}_2\equiv\int_{G({\mathbb{R}^{d}})}\varphi
	(x,y)d\bm{\nu}_1(x,y)=\int_{G({\mathbb{R}^{d}})}\varphi(x,y)d\bm{\nu
	}_2(x,y)\;,\forall\;\varphi\in C_{c}(G({\mathbb{T}^{d}})),
	\]
	where
	\[
	\mathcal{M}_{0}(\mathbb{R}^{d}):=\left\{  \bm{\nu}\in\mathcal{M}%
	_{\text{loc}}(G({\mathbb{R}^{d}}))\;|\;\int_{G(\mathbb{R}^{d})}(1\wedge
	|x-y|)\;d|\bm{\nu}|(x,y)<\infty\right\}  ,
	\]
	and $\mathcal{M}_{\text{loc}}(G({\mathbb{R}^{d}}))$ stands for the locally
	finite Radon measures on $G({\mathbb{R}^{d}})$. It is straightforward to check
	that
	\[
	\mathcal{M}_{0}(\mathbb{T}^{d})\simeq\mathcal{M}_{0}(\mathbb{R}^{d}%
	)/\bm{\sim}%
	\]
	where
	\begin{equation}
		\mathcal{M}_{0}({\mathbb{T}^{d}}):=\left\{  \bm{\nu}\in\mathcal{M}
		_{\text{loc}}(G([0,1)^{d}))\;|\;\int_{G([0,1)^{d})}d_{{\mathbb{T}^{d}}
		}(x,y)\;d|\bm{\nu}|(x,y)<\infty\right\}  , \label{M0}%
	\end{equation}
	and $\mathcal{M}_{\text{loc}}(G([0,1)^{d}))$ denotes the set of locally finite
	signed Radon measure defined in the Borelians of $G({[0,1)^{d}})$, defined in \eqref{G01}, that are
	generated by $d_{{\mathbb{T}^{d}}}$.
	
	Next, $\forall\zeta\in C({\mathbb{T}^{d}}),\; \forall[\mu]\in\mathcal{P}%
	({\mathbb{T}^{d}})$ we define $\int_{{\mathbb{T}^{d}}}\zeta\;d[\mu] :=
	\int_{{\mathbb{R}^{d}}}\zeta\; d\mu$. In particular, since there exists only
	one $\overline{\mu}\in\mathcal{P}(\mathbb{R}^{d})$ such that $\text{supp }
	\overline{\mu} \subseteq[0,1)^{d}$ e $\overline{\mu} \sim\mu$ (see \cite{CarrilloSlepcev,GT}), we actually
	have
	\[
	\int_{{\mathbb{T}^{d}}}\zeta\;d[\mu] = \int_{[0,1)^{d}} \zeta\; d\overline
	{\mu}.
	\]
	The measure $\overline{\mu}$ can be obtained by
	\[
	\overline{\mu}(B):= \sum_{k\in\mathbb{Z}^{d}}\mu(B+k)\;,\;\; \mbox{ for all
		Borelians } B\subseteq[0,1)^{d},
	\]
	which is the push-forward of $\mu$ by the map $T:{\mathbb{R}^{d}}%
	\rightarrow[0,1)$ defined by $T(x_{1},\ldots,x_{d}):=(x_{1}- \lfloor{x_{1}%
	}\rfloor,\ldots,x_{d}- \lfloor{x_{d}}\rfloor)$, where $\lfloor{.}\rfloor$ is
	the greatest integer function.
	
	Analogously, for $\varphi\in C_{c}(G({\mathbb{T}^{d}}))$and $[\bm{\nu}%
	]\in\mathcal{M}_{0}(\mathbb{T}^{d})$ we define
	\[
	\int_{G({\mathbb{T}^{d}})}\varphi\;d[\bm{\nu}]:=\int_{G({\mathbb{R}^{d}}
		)}\varphi\;d\bm{\nu}.
	\]
	In particular, the integrability condition in \eqref{M0} allows us to use the
	push-forward by $T$ and obtain that there exists $\overline{\bm{\nu}}%
	\in\mathcal{M}_{\text{loc}}(G({\mathbb{R}^{d}}))$ such that $\text{supp
	}\overline{\bm{\nu}}\subseteq G(\Q)$ and
	$\overline{\bm{\nu}}\bm{\sim}\bm{\nu}$. Therefore
	\[
	\int_{{\mathbb{T}^{d}}}\varphi\;d[\bm{\nu}]=\int_{G([0,1)^{d})}%
	\varphi\;d\overline{\bm{\nu}}.
	\]

	Let us recall the narrow convergence in $\mathcal{P}({\mathbb{T}^{d}})$ which
	will be used in the results of compactness. In fact, this topology coincides
	with the weak-$\ast$ topology induced by the dual of $C(\mathbb{T}^{d})$.
	
	\begin{definition}
		\label{Def.weakconv} We say that a sequence $(\mu_{n})_{n}\subseteq
		\mathcal{P}({\mathbb{T}^{d}})$ weakly (narrowly) converges to $\mu
		\in\mathcal{P}({\mathbb{T}^{d}})$ if
		\[
		\int_{{[0,1)^{d}}}\varphi(x)\;d\mu_{n}(x)\rightarrow\int_{{[0,1)^{d}}}
		\varphi(x)\;d\mu(x)\;,\;\;\forall\;\varphi\in C({\mathbb{T}^{d}}).
		\]
		
	\end{definition}
	
	\begin{remark}
		\label{d.Td} Note that the cube ${[0,1)^{d}}$ with the torus metric
		$d_{{\mathbb{T}^{d}}}$ (see \eqref{metric-torus-1}) is a compact metric space
		and therefore, by Prokhorov's Theorem, $\mathcal{P}({\mathbb{T}^{d}})$ is a
		weakly compact space.
	\end{remark}
	

	\section{The action functional}
	
	\quad The weak formulation given in \eqref{eq.FMPEweak} for the
	fractional porous medium equation leads us to the following function, which is
	used there as a mean for probability densities between any two different
	points. We define the $m$-mean $\theta_{m}:[0,\infty)\times\lbrack
	0,\infty)\rightarrow\lbrack0,\infty)$ by
	\begin{equation*}
		\theta_{m}(s,t):=\left\{
		\begin{array}
			[c]{cl}%
			\displaystyle\frac{s-t}{\log s-\log t}\;, & \mbox{ if }\;m=1\\
			\displaystyle\frac{m-1}{m}\frac{s^{m}-t^{m}}{s^{m-1}-t^{m-1}}\;, &
			\mbox{ else,
			}
		\end{array}
		\right.  \;\;\mbox{ for }s\neq t\;, 
	\end{equation*}
	and extended by continuity to $s=t$. The function $\theta_{m}:[0,\infty
	)^{2}\rightarrow\lbrack0,\infty)$ satisfies the following integral
	representation
	\[
	\theta_{m}(s,t)=\left\{
	\begin{array}
	[c]{ll}%
	\displaystyle\int_{0}^{1}s^{\alpha}t^{1-\alpha}\;d\alpha\;, & \mbox{ if }m=1\\
	\displaystyle\int_{0}^{1}\left(  (1-\alpha)s^{m-1}+\alpha t^{m-1}\right)
	^{\frac{1}{m-1}}\;d\alpha\;, & \mbox{ if }m\neq1
	\end{array}
	\right.
	\]
	and the following properties for $0<m\leqslant2$:
	
	\begin{enumerate}
		\item Symmetry: $\theta_{m}(s,t)=\theta_{m}(t,s)$;
		
		\item Homogeneity: $\theta_{m}(\lambda s,\lambda t)=\lambda\theta_{m}(s,t) $,
		for all $\lambda>0$;
		
		\item Concavity: $\theta_{m}$ is concave for $0<m\leqslant2$;
		
		\item Monotonicity: If $m_{1}<m_{2}$ then $\theta_{m_{1}}(s,t)<\theta_{m_{2}
		}(s,t)$ for all $s,t>0$;
		
		\item Monotonicity: If $0\leqslant s_{1}\leqslant s_{2}$ and $t\geqslant0$
		then $\theta_{m}(s_{1},t)\leqslant\theta_{m}(s_{2},t)$;
		
		\item Boundary values: $\theta_{m}(0,t)=0$ for all $t\geqslant0$ if and only
		if $m\leqslant1$. For $1<m\leqslant2$, we have $\theta_{m}(0,t)=\frac{m-1}
		{m}t$, for all $t\geqslant0$.
	\end{enumerate}
	
	Given a probability density $\rho(x)dx\in\mathcal{P}({\mathbb{T}^{d}})$, we
	use the following notation for the $m$-mean of $\rho$ between any two points
	$(x,y)\in{\mathbb{T}^{d}}\times{\mathbb{T}^{d}}$:%
	\[
	\widehat{\rho}_{m}(x,y):=\theta_{m}(\rho(x),\rho(y)).
	\]

	Let us define now the action functional that will be used in the dynamic
	definition of the non-local fractional Wasserstein distance. Just for a
	moment, let us write the torus as the cube ${[0,1)^{d}}$ with its Borelians
	defined by the metric $d_{{\mathbb{T}^{d}}}$.
	
	\begin{definition}
		Given $(\mu,\bm{\nu})\in\mathcal{P}([0,1)^{d})\times\mathcal{M}
		_{\text{loc}}(G([0,1)^{d}))$, we define $\bm{\mu}^{1},\bm{\mu}^{2}
		\in\mathcal{M}_{\text{loc}}(G([0,1)^{d}))$ by
		\[
		d\bm{\mu}^{1}(x,y)=K^{\sigma}(x-y)\;dy\;d\mu(x)\;,\;\;d\bm{\mu}
		^{2}(x,y)=K^{\sigma}(x-y)\;dx\;d\mu(y)\in\mathcal{M}_{\text{loc}}
		(G({[0,1)^{d}}))
		\]
		Now, let $\bm{\lambda}\in\mathcal{M}_{\text{loc}}(G({[0,1)^{d}}))$ be a
		nonnegative Borelian measure such that $\bm{\mu}^{1},\bm{\mu}^{2}$ and $\bm{\nu}$
		are absolutely continuous with respect to $\bm{\lambda}$ with
		$d\bm{\mu}^{1}=\rho^{1}d\bm{\lambda}$, $d\bm{\mu}^{2}=\rho
		^{2}d\bm{\lambda}$ and $d\bm{\nu}=wd\bm{\lambda}$. The measure $\lambda$ can be taken as $|\bm{\mu}^{1}|+|\bm{\mu}^{2}|+|\bm{\nu}|$, for examle. Then we define
		the action $\mathcal{A}$ by
		\begin{equation*}
			\mathcal{A}(\mu,\bm{\nu}):=\frac{1}{2}\int_{G({[0,1)^{d}})}\frac
			{w(x,y)^{2}}{\theta_{m}(\rho^{1}(x,y),\rho^{2}(x,y))}\;d\bm{\lambda}(x,y)
		\end{equation*}
		
	\end{definition}
	
	\begin{remark}
		We should mention that the action functional does depend on the order $\sigma$
		of the fractional Laplacian and also on the nonlinearity $m$, although we will
		not mention them in the notation in order to keep it clearer.
	\end{remark}
	
	This functional is well defined since the integrand is positively 1-homogeneous.
	
	Let us show that for $m=1$ the present case is the natural periodization of
	the one defined in \cite{erbar}. Given $\mu\in\mathcal{P}(\mathbb{R}^{d})$, we
	can define
	\[
	d\bm{\mu}^{1}(x,y)=\frac{C_{d,\sigma}}{|x-y|^{d+2\sigma}}dyd\mu
	(x)\in\mathcal{M}_{\text{loc}}(G({\mathbb{R}^{d}}))
	\]
	and consider the respective equivalent classes $[\mu]$ and $[\bm{\mu}%
	^{1}]$. Therefore, if $\varphi\in C_{c}(G({\mathbb{T}^{d}}))$, we can write
	for all $x\in{\mathbb{R}^{d}}$
	\begin{align*}
		\int_{{\mathbb{R}^{d}}\backslash\{x\}}\frac{\varphi(x,y)}{|x-y|^{d+2\sigma}%
		}\;dy  &  =\sum_{k\in\mathbb{Z}^{d}}\int_{([0,1)^{d}+k)\backslash\{x\}}%
		\frac{\varphi(x,y)}{|x-y|^{d+2\sigma}}\;dy=\sum_{k\in\mathbb{Z}^{d}}%
		\int_{[0,1)^{d} \backslash\{\overline{x}\}}\frac{\varphi(x,y-k)}%
		{|x-y+k|^{d+2\sigma}}\;dy\\
		&  =\int_{[0,1)^{d}\backslash\{\overline{x}\}}\varphi(x,y)\sum_{k\in
			\mathbb{Z}^{d}}\frac{1}{|x-y+k|^{d+2\sigma}}\;dy\\
		&  =\frac{1}{C_{d,\sigma}}\int_{[0,1)^{d}\backslash\{\overline{x}\}}
		\varphi(x,y)K^{\sigma}(x-y)\;dy
	\end{align*}
	and then
	\begin{align*}
		\int_{G({\mathbb{T}^{d}})}\varphi(x,y)\;d[\bm{\mu}^{1}](x,y)  &
		=\int_{G({\mathbb{R}^{d}})}\varphi(x,y)\;d\bm{\mu}^{1}(x,y)\\
		&  =\int_{{\mathbb{R}^{d}}}\int_{{\mathbb{R}^{d}}\backslash\{x\}}
		\frac{C_{d,\sigma}}{|x-y|^{d+2\sigma}}\varphi(x,y)\;dy\;d\mu(x)\\
		&  =\int_{{\mathbb{R}^{d}}}\int_{[0,1)^{d}\backslash\{\overline{x}\}}
		\varphi(x,y)K^{\sigma}(x-y)\;dy\;d\mu(x)\\
		&  =\int_{[0,1)^{d}}\int_{[0,1)^{d}\backslash\{\overline{x}\}}\varphi
		(x,y)K^{\sigma}(x-y)\;dy\;d\overline{\mu}(x)\\
		&  =\int_{G([0,1)^{d})}\varphi(x,y)K^{\sigma}(x-y)\;dy\;d\overline{\mu}(x)
	\end{align*}
	Therefore, we can write $d[\bm{\mu}^{1}](x,y)=K^{\sigma}(x-y)\;dy\;d[\mu
	](x)$. Analogously we define $d[\bm{\mu}^{2}](x,y)=K^{\sigma
	}(x-y)\;dx\;d[\mu](y)$.
	
	From now onwards we denote any equivalence class in $\mathcal{P}%
	({\mathbb{T}^{d}})$ or $\mathcal{M}_{\text{loc}}(G({\mathbb{T}^{d}}))$ by its
	representative element with support in ${[0,1)^{d}}$, and thus, we skip the
	bracket notation. Also, in order to simplify the notation, we will denote by
	$dK^{\sigma}(x,y)$ the locally finite measure $K^{\sigma}(x-y)dxdy\in
	\mathcal{M}_{\text{loc}}(G({[0,1)^{d}}))$.
	
	A key ingredient to working with the action functional and later with the
	Minimizing Movement Scheme is the lower semi-continuity property. This is a
	direct consequence of the following more general result that can be found in
	\cite{Buttazzo}.
	
	\begin{lemma}
		\label{l.s.c.Functional} Let $\Omega$ be a locally compact Polish space and
		$f:\Omega\times\mathbb{R}^{n}\rightarrow[0,+\infty]$ be a lower semicontinuous
		function such that $f(z,\cdot)$ is convex and positively 1-homogeneous for
		every $z\in\Omega$. Then the functional $F:\mathcal{M}_{\text{loc}}
		(\Omega)^{n}\rightarrow[0,+\infty]$ defined by
		\[
		F(\sigma_{1},\ldots,\sigma_{n}) = \int_{\Omega}f(z,h_{1}(z),\ldots
		,h_{n}(z))\;d\lambda(z)\;,
		\]
		where $\lambda\in\mathcal{M}_{\text{loc}}(\Omega)$ is such that $d\sigma
		_{i}=h_{i}(z)d\lambda$, for $1\leqslant i\leqslant n$, is weak* lower
		semicontinuous on $\mathcal{M}_{\text{loc}}(\Omega)$.
	\end{lemma}
	
	The results contained in the next lemma are periodic analogous to the ones
	obtained in \cite{erbar} for the equivalent problem in $\mathbb{R}^{d}$.
	
	\begin{lemma}
		\label{A.sc} The action functional $\mathcal{A}$ is convex in $\mathcal{P}%
		({\mathbb{T}^{d}})\times\mathcal{M}_{\text{loc}}(G({\mathbb{T}^{d}})).$ Also,
		it is lower semicontinuous, i.e., if $\mu_{n}\rightarrow\mu$ weakly and
		$\bm{\nu}_{n}\rightarrow\bm{\nu}$ weakly*, then
		\[
		\mathcal{A}(\mu,\bm{\nu})\leqslant\liminf_{n\rightarrow\infty}%
		\mathcal{A}(\mu_{n},\bm{\nu}_{n}).
		\]
		Furthermore, for any $d\mu=\rho(x)\;dx$ and $\bm{\nu}\in\mathcal{M}%
		_{\text{loc}}(G({\mathbb{T}^{d}}))$ such that $\mathcal{A}(\mu,\bm{\nu
		})<\infty$, there exists a function $w:G([0,1)^{d})\rightarrow\mathbb{R}$ such
		that $\bm{\nu}=w\widehat{\rho}_{m}dK^{\sigma}$ and
		\[
		\mathcal{A}(\mu,\bm{\nu})=\frac{1}{2}\int_{G([0,1)^{d})}|w(x,y)|^{2}%
		\widehat{\rho}_{m}(x,y)dK^{\sigma}(x,y).
		\]
		
	\end{lemma}
	
	\emph{Proof.} The convexity is clear from the fact that $\theta$ is concave
	and $(x,y)\mapsto x^{2}/y$ is convex and decreasing in $y$. The lower
	semicontinuity follows from Lemma \ref{l.s.c.Functional} for $\Omega
	=G({[0,1)^{d}})$, $f(z,h_{1},h_{2},h_{3}):=|h_{3}|^{2}\theta_{m}(h_{1}%
	,h_{2})^{-1}$ and noting that $\bm{\mu}^{i,n}\overset{\ast}%
	{\rightharpoonup}\bm{\mu}^{i}$ for $i=1,2$ if $\mu^{n}\rightharpoonup\mu$.
	Thus
	\[
	\mathcal{A}(\mu,\bm{\nu})=F(\bm{\mu}^{1},\bm{\mu}^{2},\bm{\nu
	})=\frac{1}{2}\int_{G({[0,1)^{d}})}f(\rho^{1},\rho^{2},w)\;d\bm{\lambda}%
	\]
	is lower semicontinuous w.r.t. the weak-weak* convergence.
	
	The second claim follows from a similar proof as the one in \cite[Lemma
	2.3]{erbar} and it is left to the reader.
	
	\fin
	
	The following lemma will be useful in order to provide uniform estimates that
	will be needed in several results in the sequel.
	
	\begin{lemma}
		\label{comp.estimate} There exists a constant $C>0$ such that
		\[
		\int_{G({[0,1)^{d}})}d_{{\mathbb{T}^{d}}}(x,y)|\bm{\nu}|(x,y)\leqslant
		C\sqrt{\mathcal{A}(\mu,\bm{\nu})},
		\]
		for all $\mu\in\mathcal{P}({\mathbb{T}^{d}})$ and $\bm{\nu}\in
		\mathcal{M}_{\text{loc}}(G({\mathbb{T}^{d}}))$. Also, for every Borelian
		$F\subseteq G({\mathbb{T}^{d}})$ such that $\delta:=\inf\{d_{{\mathbb{T}^{d}}
		}(x,y)\;:\;(x,y)\in F\}>0$, there exists a constant $C=C(\delta)$ such that
		\[
		|\bm{\nu}|(F)\leqslant C\sqrt{\mathcal{A}(\mu,\bm{\nu})}\text{.}
		\]
		
	\end{lemma}
	
	\emph{Proof. }We start by estimating
	\[
	|\bm{\nu}|(F)=\frac{1}{\delta}\int_{F}\delta\;d|\bm{\nu}%
	|\leqslant\frac{1}{\delta}\int_{F}d(x,y)\;d|\bm{\nu}|\text{.}%
	\]
	Now, let $\bm{\lambda}\in\mathcal{M}_{\text{loc}}(G({\mathbb{T}^{d}}))$
	such that $\bm{\mu}^{1},\bm{\mu}^{1},\bm{\nu}\ll\bm{\lambda}$
	with $d\bm{\mu}^{i}=\rho^{i}d\bm{\lambda}$ for $i=1,2$ and
	$d\bm{\nu}=wd\bm{\lambda}$. If $\mathcal{A}(\mu,\bm{\nu})<\infty$
	then $\bm{\lambda}\left(  \left\{  (x,y)\in G({\mathbb{T}^{d}}%
	)\;|\;\frac{w^{2}}{\theta_{m}(\rho^{1},\rho^{2})}=\infty\right\}  \right)  =0$
	and thus
	\begin{align}
		\int_{G({[0,1)^{d}})}d_{{\mathbb{T}^{d}}}(x,y)\;d|\bm{\nu}|  &
		=\int_{G({[0,1)^{d}})}d_{{\mathbb{T}^{d}}}(x,y)|w(x,y)|\;d\bm{\lambda
		}\nonumber\\
		&  =\int_{G({[0,1)^{d}})}d_{{\mathbb{T}^{d}}}(x,y)\sqrt{2\theta_{m}(\rho
			^{1},\rho^{2})}\sqrt{\frac{w^{2}}{2\theta_{m}(\rho^{1},\rho^{2})}
		}\;d\bm{\lambda}\nonumber\\
		&  \leqslant\left(  2\int_{G({[0,1)^{d}})}d_{{\mathbb{T}^{d}}}^{2}
		(x,y)\theta_{m}(\rho^{1},\rho^{2})\;d\bm{\lambda}\right)  ^{\frac{1}{2}
		}\left(  \int_{G({[0,1)^{d}})}\frac{w^{2}}{2\theta_{m}(\rho^{1},\rho^{2}
		)}\;d\bm{\lambda}\right)  ^{\frac{1}{2}}\nonumber\\
	&  =C\sqrt{\mathcal{A}(\mu,\bm{\nu})}, \label{ineq.1}%
\end{align}
where $C$ is the constant that appears in the following estimate:
\begin{align}
	\int_{G({[0,1)^{d}})}d_{{\mathbb{T}^{d}}}^{2}(x,y)\theta_{m}(\rho^{1},\rho
	^{2})\;d\bm{\lambda}  &  \leqslant\int_{G({[0,1)^{d}})}d_{{\mathbb{T}^{d}
		}}^{2}(x,y)\frac{\rho^{1}+\rho^{2}}{2}\;d\bm{\lambda}\nonumber\\
		&  \leqslant\frac{1}{2}\int_{G({[0,1)^{d}})}d_{{\mathbb{T}^{d}}}
		^{2}(x,y)K^{\sigma}(x-y)\;dyd\mu(x)\nonumber\\
		&  +\frac{1}{2}\int_{G({[0,1)^{d}})}d_{{\mathbb{T}^{d}}}^{2}(x,y)K^{\sigma
		}(x-y)\;dxd\mu(y)\nonumber\\
		&  =\int_{G({[0,1)^{d}})}d_{{\mathbb{T}^{d}}}^{2}(x,y)K^{\sigma}
		(x-y)\;dyd\mu(x)\nonumber\\
		&  \leqslant\sup_{x\in{[0,1)^{d}}}\int_{{[0,1)^{d}}}d_{{\mathbb{T}^{d}}}
		^{2}(x,y)K^{\sigma}(x-y)\;dy=:C\text{.} \label{ineq.2}%
	\end{align}
	We can see that $C$ is finite by using that $d_{{\mathbb{T}^{d}}
	}(x,y)\leqslant\min_{k\in\mathbb{Z}^{d}}\{|x-y+k|\}\leqslant1$ and
	\begin{align*}
		\int_{{[0,1)^{d}}}d_{{\mathbb{T}^{d}}}^{2}(x,y)K^{\sigma}(x-y)\;dy  &
		\leqslant C_{d,\sigma}\int_{{[0,1)^{d}}}\sum_{|k|\leqslant1}\frac
		{1}{|x-y+k|^{d+2\sigma-2}}\;dy\\
		&  +C_{d,\sigma}\int_{{[0,1)^{d}}}\sum_{|k|>1}\frac{1}{|x-y+k|^{d+2\sigma}%
		}\;dy\\
		&  \leqslant C_{d,\sigma}\int_{|x-y|\leqslant2\sqrt{d}}\frac{1}%
		{|x-y|^{d+2\sigma-2} }\;dy\\
		&  +C_{d,\sigma}\int_{|x-y|>2\sqrt{d}}\frac{1}{|x-y|^{d+2\sigma}}\;dy\\
		&  =C_{d,\sigma}\int_{|y|\leqslant2\sqrt{d}}\frac{1}{|y|^{d+2\sigma-2}
		}\;dy+C_{d,\sigma}\int_{|y|>2\sqrt{d}}\frac{1}{|y|^{d+2\sigma}}\;dy\\
		&  <\infty
	\end{align*}

	\fin
	

	\section{The non-local continuity equation in ${\mathbb{T}^{d}}$}
	
	\quad This section is devoted to present some results about the periodic continuity
	equation
	\begin{equation}
		\partial_{t}\mu_{t}+\overline{\nabla}\cdot\bm{\nu}_{t}%
		=0\;,\;\;\mbox{ in }(0,T)\times{[0,1)^{d}}, \label{PCeq}%
	\end{equation}
	where $\overline{\nabla}$ is the discret gradient defined in
	(\ref{discret-grad-1}). Here $(\mu_{t})_{t\in\lbrack0,T]}$ and $(\bm{\nu
	}_{t})_{t\in\lbrack0,T]}$ are Borel families of measures in $\mathcal{P}%
	({\mathbb{T}^{d}})$ and $\mathcal{M}_{\text{loc}}(G({\mathbb{T}^{d}}))$,
	respectively, and satisfying
	\[
	\int_{0}^{T}\int_{G({[0,1)^{d}})}d_{{\mathbb{T}^{d}}}(x,y)\;d|\bm{\nu}%
	_{t}|dt<\infty.
	\]
	We suppose that \eqref{PCeq} holds in the periodic distributional sense, i.e.,
	for all $\varphi\in C_{c}^{\infty}((0,T);C^{\infty}({\mathbb{T}^{d}}))$ we
	have that the pair $(\mu_{t},\bm{\nu}_{t})_{t\in\lbrack0,T]}$ satisfies
	\begin{equation*}
		\int_{0}^{T}\int_{{[0,1)^{d}}}\partial_{t}\varphi_{t}\;d\mu_{t}(x)\;dt+\frac
		{1}{2}\int_{0}^{T}\int_{G({[0,1)^{d}})}\overline{\nabla}\varphi_{t}%
		(x,y)\;d\bm{\nu}_{t}(x,y)\;dt=0. 
	\end{equation*}

	Let us define the class of pairs that we will use in the Benamou-Brenier formulation.
	
	\begin{definition}
		\label{definition:CE} Given $\widetilde{\mu}_{0},\widetilde{\mu}_{1}
		\in\mathcal{P}({\mathbb{T}^{d}})$, we denote $\mathcal{CE}_{T}(\widetilde{\mu
		}_{0},\widetilde{\mu}_{1})$ the set of pairs $(\mu_{t},\bm{\nu}_{t}
		)_{t\in[0,T]}$ such that
		
		\begin{itemize}
			\item[i)] $t\in[0,T]\mapsto\mu_{t}$ is weakly continuous (Definition
			\ref{Def.weakconv}) with $\mu_{0} = \widetilde{\mu}_{0}$ and $\mu_{1} =
			\widetilde{\mu}_{1}$;
			
			\item[ii)] $(\bm{\nu}_{t})_{t\in\lbrack0,T]}$ is a Borel family in
			$\mathcal{M}_{\text{loc}}(G({\mathbb{T}^{d}}))$ with $\int_{0}^{T}
			\int_{G({[0,1)^{d}})}d_{{\mathbb{T}^{d}}}(x,y)\;d|\bm{\nu}_{t}|dt<\infty$;
			
			\item[iii)] The pair $(\mu_{t},\bm{\nu}_{t})_{t\in[0,T]}$ satisfies
			$\partial_{t}\mu_{t}+\overline{\nabla}\cdot\bm{\nu}_{t}=0$ in the periodic
			distributional sense.
		\end{itemize}
	\end{definition}
	
	\begin{remark}
		Let us denote by $\mathcal{CE}_{T}$ all the pairs $(\mu_{t},\bm{\nu}
		_{t})_{t\in[0,T]}$ satisfying the items $ii)$ and $iii)$ above, with
		$t\mapsto\mu_{t}$ weakly continuous but with no fixed end points.
	\end{remark}
	
	\begin{remark}
		\label{remark.Lips} It is clear that if a pair $(\mu_{t},\bm{\nu}%
		_{t})_{t\in[0,T]}\in\mathcal{CE}_{T}(\widetilde{\mu}_{0},\widetilde{\mu}_{1})$
		and $\varphi\in C({\mathbb{T}^{d}})$ is Lipschitz, then it also satisfies
		\[
		\int\varphi(x)\;d\widetilde{\mu}_{1}(x) - \int\varphi(x)\;d\widetilde{\mu}%
		_{0}(x) = \frac{1}{2}\int_{0}^{1}\int_{G({[0,1)^{d}})} \overline{\nabla
		}\varphi(x,y)\;d\bm{\nu}_{t}dt\;.
		\]
		
	\end{remark}
	
	The condition $ii)$ ensures that the weak formulation $iii)$ is well defined
	since, for every $\varphi\in C({\mathbb{T}^{d}})$ Lipschitz and $k\in\mathbb{Z}^{d}$
	we have
	\[
	|\overline{\nabla}\varphi(x,y)|=|\varphi(x)-\varphi(y)|=|\varphi
	(x)-\varphi(y+k)|\leqslant [\varphi]_1 |x-y-k|
	\]
	where \[ [\varphi]_1 := \sup_{x\neq y}\frac{|\varphi(x)-\varphi(y)|}{|x-y|},  \]
	and therefore $|\overline{\nabla}\varphi(x,y)|\leqslant [\varphi]_1 d_{\mathbb{T}^{d}}(x,y)$.
	
	The lemma below consists in a time-rescaling property for solutions of
	(\ref{PCeq}).
	
	\begin{lemma}
		\label{rescaling} Let $\bm{t}:[0,\overline{T}]\rightarrow\lbrack0,T]$ be a
		strictly increasing absolutely continuous map with absolutely continuous
		inverse $\bm{s}:=\bm{t}^{-1}$. Then $(\mu_{t},\bm{\nu}_{t}%
		)_{t\in\lbrack0,T]}$ is a distributional solution of $\partial_{t}\mu
		_{t}+\overline{\nabla}\cdot\bm{\nu}_{t}=0$ in $(0,T)\times{[0,1)^{d}}$ in
		the periodic sense if and only if the pair $(\overline{\mu}_{s},\overline
		{\bm{\nu}}_{s})_{s\in\lbrack0,\overline{T}]}$, defined by $\overline{\mu
		}_{s}=\mu_{\bm{t}(s)}$ and $\overline{\bm{\nu}}_{s}=\bm{t}%
		^{\prime}(s)\bm{\nu}_{\bm{t}(s)}$, is a solution in $(0,\overline
		{T})\times{[0,1)^{d}}$ in the periodic sense.
	\end{lemma}
	
	\emph{Proof.} Let us assume that $\bm{s}\in C^{1}(0,\overline{T})$ with
	$\bm{s}^{\prime}>0$. For every $\overline{\varphi}\in C_{c}^{\infty
	}((0,\overline{T});C^{\infty}({\mathbb{T}^{d}}))$, we define $\varphi
	_{t}(x):=\overline{\varphi}_{\bm{s}(t)}(x)$. Therefore, we can write
	\begin{align*}
		&  \int_{0}^{T}\int_{{[0,1)^{d}}}\partial_{t}\varphi_{t}\;d\mu_{t}
		(x)\;dt+\frac{1}{2}\int_{0}^{T}\int_{G({[0,1)^{d}})}\overline{\nabla}
		\varphi_{t}(x,y)\;d\bm{\nu}_{t}(x,y)\;dt\\
		&  =\int_{0}^{T}\int_{{[0,1)^{d}}}\partial_{s}\overline{\varphi}
		_{\bm{s(t)}}(x)\;d\mu_{t}(x)\;\bm{s}^{\prime}(t)\;dt+\frac{1}{2}
		\int_{0}^{T}\int_{G({[0,1)^{d}})}\frac{\overline{\nabla}\overline{\varphi
			}_{\bm{s}(t)}(x,y)}{\bm{s}^{\prime}(t)}\;d\bm{\nu}_{t}
		(x,y)\;\bm{s}^{\prime}(t)\;dt\\
		&  =\int_{0}^{\overline{T}}\int_{{[0,1)^{d}}}\partial_{s}\overline{\varphi
		}_{s}(x)\;d\mu_{\bm{t}(s)}(x)\;ds+\frac{1}{2}\int_{0}^{\overline{T}}
		\int_{G({[0,1)^{d}})}\overline{\nabla}\overline{\varphi}_{s}(x,y)\bm{t}
		^{\prime}(s)\;d\bm{\nu}_{\bm{t}(s)}(x,y)\;ds
	\end{align*}
	And thus, the pair $(\mu_{t},\bm{\nu}_{t})_{t\in\lbrack0,T]}$ is a
	solution if and only if $(\overline{\mu}_{s},\overline{\bm{\nu}}
	_{s})_{s\in\lbrack0,\overline{T}]}$ is also a solution. \fin
	

	\section{The non-local transport distance on $\mathcal{P}({\mathbb{T}^{d}})$}
	
	\quad In this section we define the periodic non-local Wasserstein metric in the sense of
	Benamou-Brenier formula and show that in the set where it is finite, it
	defines a metric that induces a topology stronger than the narrow one (see
	Definition \ref{Def.weakconv}).
	
	\begin{definition}
		For $\mu_{0},\mu_{1}\in\mathcal{P}({\mathbb{T}^{d}})$ we define the function
		\begin{equation}
			\mathcal{W}(\mu_{0},\mu_{1})^{2}:=\inf\left\{  \left.  \int_{0}^{1}
			\mathcal{A}(\mu_{t},\bm{\nu}_{t})\;dt\;\right\vert (\mu_{t},\bm{\nu
			}_{t})_{t\in\lbrack0,1]}\in\mathcal{CE}_{1}(\mu_{0},\mu_{1})\right\}  .
			\label{eq.inf}%
		\end{equation}
		
	\end{definition}
	
	We will show later that this infimum is actually a minimum.
	
	\begin{lemma}
		\label{WT} For any $T>0$ we have
		\[
		\mathcal{W}(\mu_{0},\mu_{1})^{2}=\inf\left\{  \left.  T\int_{0}^{T}
		\mathcal{A}(\mu_{t},\bm{\nu}_{t})\;dt\;\right\vert \;(\mu_{t},\bm{\nu
		}_{t})_{t\in\lbrack0,1]}\in\mathcal{CE}_{T}(\mu_{0},\mu_{1})\right\}  .
		\]
		
	\end{lemma}
	
	\emph{Proof.} Let $\bm{s}:[0,T]\rightarrow\lbrack0,1]$ defined by
	$\bm{s}(t)=T^{-1}t$. For any $(\mu_{t},\bm{\nu}_{t})_{t\in\lbrack
		0,1]}\in\mathcal{CE}_{1}(\mu_{0},\mu_{1})$ we have, by Lemma \ref{rescaling},
	that $(\overline{\mu},\overline{\bm{\nu}})\in\mathcal{CE}_{T}(\mu_{0}
	,\mu_{1})$ where $\overline{\mu}_{t}:=\mu_{\bm{s}(t)}$ and $\overline
	{\bm{\nu}}_{t}:=T^{-1}\bm{\nu}_{\bm{s}(t)}$. Now, since
	$\mathcal{A}$ is $2-$homogeneous in the variable $\bm{\nu}$, we have
	\[
	\int_{0}^{T}\mathcal{A}(\overline{\mu}_{t},\overline{\bm{\nu}}
	_{t})\;dt=\frac{1}{T^{2}}\int_{0}^{T}\mathcal{A}(\mu_{\bm{s}
		(t)},\bm{\nu}_{\bm{s}(t)})\;dt=\frac{1}{T}\int_{0}^{1}\mathcal{A}
	(\mu_{t},\bm{\nu}_{t})\;dt.
	\]
	Using the inverse $\bm{s}^{-1}$ we can also see that $\mathcal{CE}_{T}
	(\mu_{0},\mu_{1})$ is mapped onto $\mathcal{CE}_{1}(\mu_{0},\mu_{1})$ with the
	same relation as above. This concludes the proof. \fin
	
	\begin{lemma}
		\label{WT2} For any $T>0$ and $\mu_{0},\mu_{1}\in\mathcal{P}({\mathbb{T}^{d}
		})$, we have
		\begin{equation}
			\mathcal{W}(\mu_{0},\mu_{1})=\inf\left\{  \left.  \int_{0}^{T}\sqrt
			{\mathcal{A}(\mu_{t},\bm{\nu}_{t})}\;dt\;\right\vert \;(\mu_t,\bm{\nu
			}_t)_{t\in[0,T]}\in\mathcal{CE}_{T}(\mu_{0},\mu_{1})\right\}  .
			\label{Wbar}%
		\end{equation}
		
	\end{lemma}
	
	\emph{Proof.} Let us denote by $\widehat{\mathcal{W}}(\mu_{0},\mu_{1})$ the
	right-hand side of \eqref{Wbar}. For every $(\mu_t,\bm{\nu}_t)_{t\in
		\lbrack0,1]}\in\mathcal{CE}_{1}(\mu_{0},\mu_{1})$ we have by Cauchy-Schwarz
	inequality that%
	
	\[
	\int_{0}^{T}\sqrt{\mathcal{A}(\mu_{t},\bm{\nu}_{t})}\;dt\leqslant\sqrt
	{T}\left(  \int_{0}^{T}\mathcal{A}(\mu_{t},\bm{\nu}_{t})\;dt\right)
	^{\frac{1}{2}}\;,
	\]
	and thus, by Lemma \ref{WT}, $\widehat{\mathcal{W}}(\mu_{0},\mu
	_{1})\leqslant\mathcal{W}(\mu_{0},\mu_{1})$. For the opposite inequality we
	argue as in \cite{DNS}. Let us suppose that $\mathcal{W}(\mu_{0},\mu_{1})<\infty$ and
	let $(\mu_t,\bm{\nu}_t)_{t\in[0,1]}\in\mathcal{CE}_{1}(\mu_{0},\mu_{1})$
	be such that
	\begin{equation}\label{A.finite}
	\int_{0}^{T}\mathcal{A}(\mu_{t},\bm{\nu}_{t})\;dt<\infty\;.
	\end{equation}
	Then, for every $\varepsilon>0$ we define the function $\bm{s}
	_{\varepsilon}$ on $(0,T)$ by
	\[
	\bm{s}_{\varepsilon}(t):=\int_{0}^{t}\sqrt{\varepsilon+\mathcal{A}(\mu
		_{r},\bm{\nu}_{r})}\;dr\;.
	\]
	We have that $\bm{s}_{\varepsilon}$ is strictly increasing and its
	derivative exists a.e. with $\bm{s}_{\varepsilon}^{\prime}\geqslant
	\sqrt{\varepsilon}$. So its inverse $\bm{t}_{\varepsilon}:=\bm{s}
	_{\varepsilon}^{-1}$ is well defined and satisfies the hypothesis from Lemma
	\ref{rescaling}. Therefore, we know that $(\widehat{\mu}_{s},\widehat
	{\bm{\nu}}_{s})_{s\in\lbrack0,\bm{s}_{\varepsilon}(T)]}\in
	\mathcal{CE}_{\bm{s}_{\varepsilon}(T)}(\mu_{0},\mu_{1})$ where
	$\widehat{\mu}_{s}:=\mu_{\bm{t}_{\varepsilon}(s)}$ and $\widehat
	{\bm{\nu}}_{s}:=\bm{t}_{\varepsilon}^{\prime}(s)\bm{\nu
	}_{\bm{t}_{\varepsilon}(s)}$. Using Lemma \ref{WT}, we have
	\begin{align*}
		\mathcal{W}(\mu_{0},\mu_{1})^{2}  &  \leqslant\bm{s}_{\varepsilon}
		(T)\int_{0}^{\bm{s}_{\varepsilon}(T)}\mathcal{A}(\widehat{\mu}
		_{s},\widehat{\bm{\nu}}_{s})\;ds\\
		&  =\bm{s}_{\varepsilon}(T)\int_{0}^{\bm{s}_{\varepsilon}
			(T)}\bm{t}_{\varepsilon}^{\prime}(s)^{2}\mathcal{A}(\mu_{\bm{t}
			_{\varepsilon}(s)},\bm{\nu}_{\bm{t}_{\varepsilon}(s)})\;ds\\
		&  =\bm{s}_{\varepsilon}(T)\int_{0}^{T}(\bm{t}_{\varepsilon}^{\prime
		}\circ\bm{s}_{\varepsilon}(t))^{2}\mathcal{A}(\mu_{t},\bm{\nu}
		_{t})\bm{s}_{\varepsilon}^{\prime}(t)\;dt\\
		&  =\bm{s}_{\varepsilon}(T)\int_{0}^{T}\frac{\mathcal{A}(\mu
			_{t},\bm{\nu}_{t})}{\varepsilon+\mathcal{A}(\mu_{t},\bm{\nu}_{t}
			)}\sqrt{\varepsilon+\mathcal{A}(\mu_{t},\bm{\nu}_{t})}\;dt\\
		&  \leqslant\bm{s}_{\varepsilon}(T)^{2}%
	\end{align*}
	This holds for all $\varepsilon>0$ and therefore
	\[
	\mathcal{W}(\mu_{0},\mu_{1})^{2}\leqslant\lim_{\varepsilon\rightarrow0^{+}
	}\bm{s}_{\varepsilon}(T)^{2}=\left(  \int_{0}^{T}\mathcal{A}(\mu
	_{t},\bm{\nu}_{t})\;dt\right)  ^{2}\;.
	\]
	Since we have this inequality for all pairs $(\mu_t,\bm{\nu}_t)_{t\in\lbrack0,T]}$ satisfying \eqref{A.finite}, we can take the infimum over $\mathcal{CE}_T(\mu_0,\mu_1)$ and obtain $\mathcal{W}(\mu_{0},\mu_{1})\leqslant\widehat{\mathcal{W}}(\mu_{0},\mu_{1})$. \fin
	
	\begin{proposition}
		\label{Wmetric} Given $\mu_{\ast}\in\mathcal{P}({\mathbb{T}^{d}})$, the
		function $\mathcal{W}$ is a metric in $\mathcal{P}_{\mu_{\ast}}:=\{\mu
		\in\mathcal{P}({\mathbb{T}^{d}})\;|\;\mathcal{W}(\mu_{\ast},\mu)<\infty\}.$
	\end{proposition}
	
	\emph{Proof.} \textit{Symmetry:} Let $\mu_{0},\mu_{1}\in\mathcal{P}_{\mu
		_{\ast}}$ and $(\mu_{t},\bm{\nu}_{t})_{t\in\lbrack0,1]}\in\mathcal{CE}_{1}(\mu
	_{0},\mu_{1})$ and consider the pair $(\widetilde{\mu}_{t},\widetilde{\bm{\nu}
	}_{t}):=(\mu_{1-t},-\bm{\nu}_{1-t})$. It is straightforward to check that $(\widetilde{\mu}_{t},\widetilde{\bm{\nu}}_{t})\in\mathcal{CE}_{1}(\mu_{1},\mu_{0})$ (see Lemma \ref{rescaling}). Since
	\[
	\int_{0}^{1}\mathcal{A}(\widetilde{\mu}_{t},\widetilde{\nu}_{t})\ dt=\int
	_{0}^{1}\mathcal{A}(\mu_{1-t},-\nu_{1-t})\ dt=\int_{0}^{1}\mathcal{A}(\mu
	_{t},\nu_{t})\ dt,
	\]
	we have that
	\[
	\mathcal{W}(\mu_{1},\mu_{0})\leqslant\mathcal{W}(\mu_{0},\mu_{1}).
	\]
	Analogously, it follows the opposite inequality.
	
	\textit{Triangular Inequality:} Assume that $\mu_{2}\in\mathcal{P}_{\mu_{\ast
		}}$. Choose $(\bar{\mu}_{t},\bar{\bm{\nu}}_{t})\in\mathcal{CE}_{1}(\mu_{0},\mu
		_{1})$ and $(\bar{\bar{\mu}}_{t},\bar{\bar{\bm{\nu}}}_{t})\in\mathcal{CE}_{1}
		(\mu_{1},\mu_{2})$. Define the pair
		\[
		(\widetilde{\mu}_{t},\widetilde{\bm{\nu}}_{t})=\left\{
		\begin{array}
		[c]{lcc}%
		(\bar{\mu}_{2t},\bar{\bm{\nu}}_{2t}) & \text{ if } & t\in\lbrack0,1/2)\\
		(\bar{\bar{\mu}}_{2t-1},\bar{\bar{\bm{\nu}}}_{2t-1}) & \text{ if } & t\in
		\lbrack1/2,1]
		\end{array}
		\right.  .
		\]
		Let us to check that $(\widetilde{\mu}_{t},\widetilde{\bm{\nu}}_{t})$ satisfies the non-local continuity equation in the periodic distributional sense according to Definition
		\ref{definition:CE}(iii). Let $\varphi\in C_{c}^{\infty}((0,1);C^{\infty
		}(\mathbb{T}^{d}))$, then
		\begin{align*}
			&  \ \int_{0}^{1}\int_{[0,1)^{d}}\partial_{t}\varphi_{t}\ d\widetilde{\mu}_{t}
			(x)\ dt+\frac{1}{2}\int_{0}^{1}\int_{[0,1)^{d}}\bar{\nabla}\varphi_{t}
			\ d\widetilde{\bm{\nu}}_{t}(x,y)\ dt\\
			&  =\frac{1}{2}\int_{0}^{1}\int_{[0,1)^{d}}\partial_{t}\varphi_{t/2}
			\ d\bar{\mu}_{t}(x)\ dt+\frac{1}{4}\int_{0}^{1}\int_{[0,1)^{d}}\bar{\nabla
			}\varphi_{t/2}\ d\bar{\bm{\nu}}_{t}(x,y)\ dt\\
			&  \ +\frac{1}{2}\int_{0}^{1}\int_{[0,1)^{d}}\partial_{t}\varphi
			_{(t+1)/2}\ d\bar{\bar{\mu}}_{t}(x)\ dt+\frac{1}{4}\int_{0}^{1}\int
			_{[0,1)^{d}}\bar{\nabla}\varphi_{(t+1)/2}\ d\bar{\bar{\bm{\nu}}}_{t}(x,y)\ dt\\
			&  =\frac{1}{2}\left(  \int_{[0,1)^{d}}\partial_{t}\varphi_{1/2}\ d\bar{\mu
			}_{1}-\int_{[0,1)^{d}}\partial_{t}\varphi_{1/2}\ d\bar{\bar{\mu}}_{0}\right)
			=0,
		\end{align*}
		because $\bar{\mu}_{1}=\bar{\bar{\mu}}_{0}=\mu_{1}$. Thus, $(\mu_{t},\bm{\nu}
		_{t})\in\mathcal{CE}_{1}(\mu_{0},\mu_{2})$ and
		\[
		\mathcal{W}(\mu_{0},\mu_{2})\leqslant\int_{0}^{1}\mathcal{A}(\bar{\mu}
		_{t},\bar{\bm{\nu}}_{t})\ dt+\int_{0}^{1}\mathcal{A}(\bar{\bar{\mu}}_{t},\bar
		{\bar{\bm{\nu}}}_{t})\ dt.
		\]
		Taking the infimum, it follows the triangular inequality.
		
		Now, we assume that $\mathcal{W}(\mu_{0},\mu_{1})=0$. By definition, we can
		get a sequence $(\mu_{t}^{k},\bm{\nu}_{t}^{k})_{k\in\mathbb{N}}\in\mathcal{CE}%
		_{1}(\mu_{0},\mu_{1})$ such that
		\[
		\lim_{k\rightarrow0}\int_{0}^{1}\mathcal{A}(\mu_{t}^{k},\bm{\nu}_{t}^{k})\ dt=0.
		\]
		Now, for any $\varphi\in C({\mathbb{T}^{d}})$ Lipschitz, we have from Remark
		\ref{remark.Lips} and from Lemma \ref{comp.estimate} that the following
		estimate holds for all $k\in\mathbb{N}$:
		\begin{align*}
			\left\vert \int\varphi\;d\mu_{1}-\int\varphi\;d\mu_{0}\right\vert  &
			\leqslant\frac{1}{2}\int_{0}^{1}\int_{G({[0,1)^{d}})}|\overline{\nabla}%
			\varphi|\;d|\bm{\nu}_{t}^{k}|dt\\
			&  \leqslant\frac{\lbrack\varphi]_{\emph{Lip}}}{2}\int_{0}^{1}\int
			_{G({[0,1)^{d}})}d_{{\mathbb{T}^{d}}}(x,y)\;d|\bm{\nu}_{t}^{k}|dt\\
			&  \leqslant C\frac{[\varphi]_{\emph{Lip}}}{2}\int_{0}^{1}\sqrt{\mathcal{A}%
				(\mu_{t}^{k},\nu_{t}^{k})}\;dt\\
			&  \leqslant C\frac{[\varphi]_{\emph{Lip}}}{2}\left(  \int_{0}^{1}%
			\mathcal{A}(\mu_{t}^{k},\nu_{t}^{k})\;dt\right)  ^{1/2}\;,
		\end{align*}
		where $[\varphi]_{\emph{Lip}}$ is the Lipschitz constant of $\varphi$.
		Therefore, taking the limit as $k\rightarrow\infty$, we obtain
		\[
		\int\varphi\;d\widetilde{\mu}_{0}=\int\varphi\;d\widetilde{\mu}_{1}%
		\]
		for all Lipschitz $\varphi\in C({\mathbb{T}^{d}})$. By an approximation
		argument, we conclude that $\mu_{0}=\mu_{1}$. \fin
		
		For the next result, let us recall the following fact about the
		Kantorovich-Rubinstein metric (see \cite[Theorem 1.14]{villani}):
		
		\begin{lemma}
			\label{Kantorovich} Let $(X,d)$ be a Polish metric space. The
			Kantorovich-Rubinstein distance on $\mathcal{P}(X)$ is defined by
			\[
			W_{d}(\mu_{0},\mu_{1}) := \inf_{\pi\in\Gamma(\mu_{0},\mu_{1})} \int_{X\times
				X}d(x,y)\;d\pi(x,y)
			\]
			where $\Gamma(\mu_{0},\mu_{1})$ is the set of probability measures $\pi$ on $X\times X$ with
			marginals $\mu_{0}$ and $\mu_{1}$. Then we have:
			
			\begin{enumerate}
				[i)]
				
				\item (Kantorovich-Rubinstein Theorem) Let $\emph{Lip}_{1}(X) := \{\varphi\in
				C(X)\;|\; [\varphi]_{\emph{Lip}}< 1 \}$ where
				\[
				[\varphi]_{\emph{Lip}}:= \sup_{x\neq y}\frac{|\varphi(x)-\varphi(y)|}{d(x,y)}%
				\]
				Then
				\[
				W_{d}(\mu_{0},\mu_{1}) = \sup_{\varphi\in\emph{Lip}_{1}(x)} |\int_{X}
				\varphi\;d\mu_{0} - \int_{X} \varphi\;d\mu_{1}|
				\]

				\item $W_{d}$ is lower semicontinuous w.r.t. weak (narrow) convergence in
				$\mathcal{P}(X)$.
			\end{enumerate}
		\end{lemma}
		
		\begin{proposition}
			[Compactness of solutions]\label{compactness} Let $\{(\mu_{t}^{n},\bm{\nu
			}_{t}^{n})_{t\in\lbrack0,T]}\}_{n\in\mathbb{N}}$ be a sequence in
			$\mathcal{CE}_{T}$ such that
			\[
			S:=\sup_{n\in\mathbb{N}}\int_{0}^{T}\mathcal{A}(\mu_{t}^{n},\bm{\nu}
			_{t}^{n})\;dt<\infty.
			\]
			Then there exists a subsequence $\{(\mu_{t}^{n_{k}},\bm{\nu}_{t}^{n_{k}
			})_{t\in\lbrack0,T]}\}_{k\in\mathbb{N}}$ and a pair $(\mu_{t},\bm{\nu}
			_{t})_{t\in\lbrack0,T]}\in\mathcal{CE}_{T}$ such that
			\begin{align*}
				\mu_{t}^{n_{k}}  &  \rightharpoonup\mu_{t}\;\mbox{ weakly in }\mathcal{P}
				({\mathbb{T}^{d}})\;\mbox{ for all }t\in\lbrack0,T]\\
				\bm{\nu}^{n_{k}}  &  \overset{\ast}{\rightharpoonup}\bm{\nu}\;\mbox{
					weakly-$*$ in }\mathcal{M}_{\text{loc}}(G({\mathbb{T}^{d}})\times
				\lbrack0,T])
			\end{align*}
			Furthermore, for any convergent sequence of couples $\{(\mu_{t}^{n}
			,\bm{\nu}_{t}^{n})_{t\in\lbrack0,T]}\}_{n\in\mathbb{N}}$ in the above
			sense, we have the following lower semicontinuity formula
			\begin{equation}
				\int_{0}^{T}\mathcal{A}(\mu_{t},\bm{\nu}_{t})\;dt\leqslant\liminf
				_{n\rightarrow\infty}\int_{0}^{T}\mathcal{A}(\mu_{t}^{n},\bm{\nu}_{t}
				^{n})\;dt. \label{int.A.lsc}%
			\end{equation}
			
		\end{proposition}
		
		\emph{Proof.} Let us first define for each $n\in\mathbb{N}$ the measures
		$d\bm{\nu}_{n}(x,y,t):=d\bm{\nu}^{n}(x,y)dt$. So, for every compact
		subset $K\subset G({\mathbb{T}^{d}})$ and every Borelian $B\subset\lbrack
		0,T]$, we can use Lemma \ref{comp.estimate} to obtain a constant $C(K)$
		depending on $K$ and write
		\begin{align}
			\sup_{n\in\mathbb{N}}|\bm{\nu}_{n}|(K\times B)  &  =\sup_{n\in\mathbb{N}
			}\int_{B}|\bm{\nu}_{t}^{n}|(K)\;dt\leqslant C(K)\sup_{n\in\mathbb{N}}
			\int_{B}\sqrt{\mathcal{A}(\mu_{t}^{n},\bm{\nu}_{t}^{n})}
			\;dt\label{disintegration}\\
			&  \leqslant C(K)|B|^{1/2}\sup_{n\in\mathbb{N}}\left(  \int_{0}^{T}
			\mathcal{A}(\mu_{t}^{n},\bm{\nu}_{t}^{n})\;dt\right)  ^{1/2}\nonumber\\
			&  \leqslant C(K)|B|^{1/2}S^{1/2}.\nonumber
		\end{align}
		Therefore, the measures $\bm{\nu}_{n}$ have total variation uniformly
		bounded on each compact subset of $G({\mathbb{T}^{d}})\times\lbrack0,T]$,
		which assures that we can extract a subsequence $\bm{\nu}_{n_k}$, such that $\bm{\nu}_{n_k}\overset{\ast}{\rightharpoonup
		}\bm{\nu}$ for some $\bm{\nu}\in\mathcal{M}_{\text{loc}}
		(G({\mathbb{T}^{d}})\times\lbrack0,T])$ as $k\to\infty$. Furthermore, the estimate
		\eqref{disintegration} leads to
		\[
		\bm{\nu}(K\times B)\leqslant\int_{B}m_{K}(t)\;dt\;,\;\;\mbox{ for }m_{K}
		(t):=\sup_{n\in\mathbb{N}}|\bm{\nu}_{t}^{n}|(K)\in L^{1}(0,T),
		\]
		so, by the Disintegration Theorem (see \cite[Theorem 2.28]%
		{ambrosio-fusco-pallara}), there exists a Borel family $(\bm{\nu}
		_{t})_{t\in\lbrack0,T]}\subseteq\mathcal{M}_{\text{loc}}(G({\mathbb{T}^{d}}))$
		such that $\bm{\nu}$ can be disintegrated w.r.t. the Lebesgue measure on
		$[0,T]$ as $d\bm{\nu}(x,y,t)=d\bm{\nu}_{t}(x,y)dt$.
		
		Let us show now that the subsequence $\{ \mu^{n_k}_{t} \}_{n}$ also admits
		subsequence which is weakly convergent for all $t\in[0,T]$. For any $k\in \mathbb{N}$,
		$0\leqslant s\leqslant t\leqslant T$ and $\varphi\in C({\mathbb{T}^{d}})$
		Lipschitz, we can use the same argument as in the proof of Proposition
		\ref{Wmetric} to obtain
		\begin{align*}
			\left|  \int\varphi\;d\mu^{n_k}_{t} - \int\varphi\;d\mu^{n_k}_{s} \right|   &
			\leqslant C\frac{[\varphi]_{\emph{Lip}}}{2} \int_{s}^{t} \sqrt{\mathcal{A}%
				(\mu_{t}^{n_k},\bm{\nu}_{t}^{n_k})}\;dt\\
			&  \leqslant C\frac{[\varphi]_{\emph{Lip}}}{2} \left(  \int_{0}^{1}
			\mathcal{A}(\mu_{t}^{n_k},\bm{\nu}_{t}^{n_k})\;dt\right)  ^{1/2}\sqrt{t-s}\\
			&  \leqslant C\frac{[\varphi]_{\emph{Lip}}}{2} S^{1/2}\sqrt{t-s}%
		\end{align*}
		Therefore, by the Kantorovich-Rubinstein Theorem (see Lemma \ref{Kantorovich}
		$(i)$), we have that
		\begin{equation}
			\label{limsupWd}\limsup_{k\to\infty}W_{d}(\mu_{s}^{n_k},\mu_{t}^{n_k}) \leqslant
			CS^{1/2}|t-s|^{1/2}%
		\end{equation}
		Now, by the compactness of $\PT$ w.r.t. the weak convergence (Remark \ref{d.Td})
		and the item $ii)$ from Lemma \ref{Kantorovich}, we have that the weak
		topology on $\mathcal{P}({\mathbb{T}^{d}})$ is compatible with $W_{d}$.
		This fact together with \eqref{limsupWd} allow us to apply \cite[Proposition 3.3.1]%
		{ambrosiogiglisavare} to the sequence $\{\mu_{t}^{n_k}\}_{k}$ and obtain that
		there exists a subsequence, still denoted by $\{\mu_{t}^{n_{k}}\}_{k}$, and a $W_{d}$-continuous
		curve $\mu:[0,T]\rightarrow\mathcal{P}({\mathbb{T}^{d}})$ such that $\mu
		_{t}^{n_{k}}\rightharpoonup\mu_{t}$ weakly for all $t\in\lbrack0,T]$ as $k\to\infty$.
		
		In order to see that the pair $(\mu_{t},\bm{\nu} _{t})_{t\in\lbrack0,T]}$
		belongs to $\mathcal{CE}_{T}$, we first note that by Lemma \ref{comp.estimate}
		and the hypothesis above, we have the following estimate
		\begin{align}
			\int_{0}^{t}\int_{G({[0,1)^{d}})}d_{{\mathbb{T}^{d}}}(x,y)\;d|\bm{\nu}
			_{t}|dt  &  \leqslant\sup_{n\in\mathbb{N}}\int_{0}^{t}\int_{G({[0,1)^{d}}
				)}d_{{\mathbb{T}^{d}}}(x,y)\;d|\bm{\nu}_{t}^{n}|dt\label{L1}\\
			&  \leqslant C\sup_{n\in\mathbb{N}}\int_{0}^{T}\sqrt{\mathcal{A}(\mu_{t}
				^{n},\bm{\nu}_{t}^{n})}\;dt\nonumber\\
			&  \leqslant CT^{1/2}\sup_{n\in\mathbb{N}}\left(  \int_{0}^{T}\mathcal{A}
			(\mu_{t}^{n},\bm{\nu}_{t}^{n})\;dt\right)  ^{1/2}\nonumber\\
			&  \leqslant CT^{1/2}S^{1/2}<\infty.\nonumber
		\end{align}
		Also, since $\mu:[0,T]\rightarrow\mathcal{P}({\mathbb{T}^{d}})$ is $W_{d}%
		$-continuous, it implies that it is weakly continuous as well.
		
		Let us now prove that $(\mu_{t},\bm{\nu}_{t})_{t\in\lbrack0,T]}$ satisfies
		the continuity equation in the periodic distributional sense. Let $\varphi\in
		C_{c}^{\infty}(0,T;C^{\infty}({\mathbb{T}^{d}}))$. Then, since $\overline
		{\nabla}\varphi\not \in C_{c}(G({[0,1)^{d}})\times\lbrack0,T])$, we use an
		argument by approximation. Let $\eta\in C_{c}^{\infty}(-1,1)$ be such that
		$0\leqslant\eta\leqslant1$ and $\eta(s)=1$ for $|s|\leqslant1/2$ and, for all
		$\varepsilon>0$, let us define $\eta_{\varepsilon}:G({[0,1)^{d}}%
		)\rightarrow\lbrack0,1]$ by $\eta_{\varepsilon}(x,y)=\eta(d_{{\mathbb{T}^{d}}%
		}(x,y)/\varepsilon)$. Thus, the function $(1-\eta_{\varepsilon})\overline
		{\nabla}\varphi$ does belong to $C_{c}(G({[0,1)^{d}})\times\lbrack0,T])$ and
		we have
		\begin{equation}
			\lim_{k\rightarrow\infty}\int_{0}^{T}\int_{G({[0,1)^{d}})}(1-\eta
			_{\varepsilon})\overline{\nabla}\varphi\;d\bm{\nu}_{t}^{n_{k}}dt=\int
			_{0}^{T}\int_{G({[0,1)^{d}})}(1-\eta_{\varepsilon})\overline{\nabla}%
			\varphi\;d\bm{\nu}_{t}dt\;,\;\;\mbox{ for all }\;\varepsilon>0.
			\label{lim1}%
		\end{equation}
		Let us define now the set $D_{\varepsilon}:=\{(x,y)\in G({[0,1)^{d}%
		})\;|\;d_{{\mathbb{T}^{d}}}(x,y)\leqslant\varepsilon\}$. Since $\text{supp
	}\eta_{\varepsilon}\subseteq D_{\varepsilon}$, we obtain
	\begin{align*}
		\left\vert \int_{0}^{T}\int_{G({[0,1)^{d}})}\eta_{\varepsilon}\overline
		{\nabla}\varphi_{t}\;d\bm{\nu}_{t}^{n_{k}}dt\right\vert  &  \leqslant
		\int_{0}^{T}\left\Vert \nabla\varphi_{t}\right\Vert _{\infty}\int
		_{D_{\varepsilon}}d_\Td(x,y)\;d|\bm{\nu}_{t}^{n_{k}}|dt\\
		&  \leqslant\sup_{t\in\lbrack0,T]}\left\Vert \nabla\varphi_{t}\right\Vert
		_{\infty}\int_{0}^{T}\int_{D_{\varepsilon}}d_{{\mathbb{T}^{d}}}%
		(x,y)\;d|\bm{\nu}_{t}^{n_{k}}|dt.
	\end{align*}
	Arguing as in \eqref{ineq.1} with $D_{\varepsilon}$ in the place of
	$G({[0,1)^{d}})$, we obtain with the inequality \eqref{ineq.2} that
	\begin{align}
		\int_{0}^{T}\int_{D_{\varepsilon}}d_{{\mathbb{T}^{d}}}(x,y)\;d|\bm{\nu
		}_{t}^{n_{k}}|dt  &  \leqslant\sqrt{2}\int_{0}^{T}\sqrt{\mathcal{A}(\mu
		_{t}^{n_{k}},\bm{\nu}_{t}^{n_{k}})}\left(  \int_{D_{\varepsilon}%
	}d_{{\mathbb{T}^{d}}}^{2}(x,y)K^{\sigma}(x-y)\;dyd\mu_{t}^{n_{k}}(x)\right)
	^{1/2}\;dt\nonumber\\
	&  \leqslant\sqrt{2}\left(  \int_{0}^{T}\mathcal{A}(\mu_{t}^{n_{k}%
	},\bm{\nu}_{t}^{n_{k}})dt\right)  ^{1/2}\nonumber\\
	&  \times\left(  \int_{0}^{T}\int_{D_{\varepsilon}}d_{{\mathbb{T}^{d}}}%
	^{2}(x,y)K^{\sigma}(x-y)\;dyd\mu_{t}^{n_{k}}(x)dt\right)  ^{1/2}\nonumber\\
	&  \leqslant\sqrt{2}S^{1/2}T^{1/2}\left(  \sup_{x\in{[0,1)^{d}}}%
	\int_{\{y|d(x,y)\leqslant\varepsilon\}}d_{{\mathbb{T}^{d}}}^{2}(x,y)K^{\sigma
	}(x-y)\;dy\right)  ^{1/2}. \label{ineq.4}%
\end{align}
Since the right-hand side of \eqref{ineq.4} goes to zero as $\varepsilon
\rightarrow0$, we have that
\begin{equation}
	\lim_{\varepsilon\rightarrow0}\sup_{k\in\mathbb{N}}\left\vert \int_{0}^{T}%
	\int_{G({[0,1)^{d}})}\eta_{\varepsilon}\overline{\nabla}\varphi_{t}%
	\;d\bm{\nu}_{t}^{n_{k}}dt\right\vert =0. \label{lim2}%
\end{equation}
Now, since $(1-\eta_{\varepsilon})\overline{\nabla}\varphi_{t}(x,y)\rightarrow
\overline{\nabla}\varphi_{t}(x,y)$, for all $(x,y)\in G({[0,1)^{d}})$ and
$t\in\lbrack0,T]$, and
\[
|(1-\eta_{\varepsilon})\overline{\nabla}\varphi_{t}(x,y)|\leqslant\sup
_{t\in\lbrack0,T]}\left\Vert \nabla\varphi_{t}\right\Vert _{\infty
}d_{{\mathbb{T}^{d}}}(x,y)\in L^{1}(d\nu_{t}dt)\;,\;\mbox{ by }\eqref{L1}\;,
\]
we can use the Dominated Convergence Theorem to obtain
\begin{equation}
	\lim_{\varepsilon\rightarrow0}\int_{0}^{T}\int_{G({[0,1)^{d}})}(1-\eta
	_{\varepsilon})\overline{\nabla}\varphi_{t}(x,y)\;d\bm{\nu}_{t}dt=\int
	_{0}^{T}\int_{G({[0,1)^{d}})}\overline{\nabla}\varphi_{t}(x,y)\;d\bm{\nu
	}_{t}dt. \label{lim3}%
\end{equation}
Therefore, from \eqref{lim1}, \eqref{lim2} and \eqref{lim3} we conclude that
\[
\lim_{k\rightarrow\infty}\int_{0}^{T}\int_{G({[0,1)^{d}})}\overline{\nabla
}\varphi_{t}(x,y)\;d\bm{\nu}_{t}^{n_{k}}dt=\int_{0}^{T}\int_{G({[0,1)^{d}%
})}\overline{\nabla}\varphi_{t}(x,y)\;d\bm{\nu}_{t}dt.
\]
Finally, $\mu_{t}^{n_{k}}\rightharpoonup\mu_{t}$ weakly for all $t\in
\lbrack0,T]$ implies that
\[
\lim_{k\rightarrow\infty}\int_{0}^{T}\int_{{[0,1)^{d}}}\partial_{t}\varphi
_{t}(x)\;d\mu_{t}^{n_{k}}(x)dt=\int_{0}^{T}\int_{{[0,1)^{d}}}\partial
_{t}\varphi_{t}(x)\;d\mu_{t}(x)dt
\]
and this shows that $(\mu_{t},\bm{\nu}_{t})_{t\in\lbrack0,T]}$ satisfies
the continuity equation in the periodic distributional sense.

The lower semicontinuity property follows from Lemma \ref{l.s.c.Functional} by
considering $\Omega=G({\mathbb{T}^{d}})\times\lbrack0,T]$, the function
$f(z,h_{1},h_{2},h_{3}):=|h_{3}|^{2}\theta_{m}(h_{1},h_{2})^{-1}$ and
functional $F:\left(  \mathcal{M}_{\text{loc}}(G({\mathbb{T}^{d}}%
)\times\lbrack0,T])\right)  ^{3}\rightarrow\lbrack0,\infty]$ defined by
\[
(\mu_{t},\bm{\nu}_{t})_{t\in\lbrack0,T]}\mapsto F(\bm{\mu}_{t}%
^{1},\bm{\mu}_{t}^{2},\bm{\nu}_{t}):=\int_{0}^{T}\int_{G({[0,1)^{d}}
	)}f(\rho_{t}^{1},\rho_{t}^{2},w_{t})\;d\bm{\lambda}_{t}dt.
\]
Therefore, we have that $F$ is weakly* lower semicontinuous, which gives
\eqref{int.A.lsc} because $\mu_{t}^{n}\rightharpoonup\mu_{t}$ implies that
$\bm{\mu}_{t}^{i,n}dt\overset{\ast}{\rightharpoonup}\bm{\mu}_{t}%
^{i}dt$ weakly* in $\mathcal{M}_{\text{loc}}(G({\mathbb{T}^{d}})\times
\lbrack0,T])$. \fin

\begin{proposition}
	\label{geodesic} Let $\mu_{0},\mu_{1}\in\mathcal{P}({\mathbb{T}^{d}})$ be such
	that $\mathcal{W}(\mu_{0},\mu_{1})<\infty$. Then there exists a pair $(\mu
	_{t},\bm{\nu}_{t})_{t\in\lbrack0,1]}\subset\mathcal{CE}_{1}(\mu_{0}
	,\mu_{1})$ which attains the infimum in \eqref{eq.inf}. Furthermore, any such
	curve satisfies
	\begin{equation}
		\mathcal{A}(\mu_{t},\bm{\nu}_{t})=\mathcal{W}(\mu_{0},\mu_{1}
		)^{2}\;,\;\;\mbox{
			for almost every }\;\;t\in\lbrack0,1]. \label{W=A}%
	\end{equation}
	Also, any curve $(\mu_{t})_{t\in\lbrack0,1]}$ satisfying \eqref{W=A} is a
	constant speed geodesic for $\mathcal{W}$, i.e.
	\[
	\mathcal{W}(\mu_{s},\mu_{t})=(t-s)\mathcal{W}(\mu_{0},\mu_{1}
	)\;,\;\;\mbox{ for
		all }\;\;0\leqslant s\leqslant t\leqslant1.
	\]
	
\end{proposition}

\emph{Proof.} If $\mu_{0},\mu_{1}\in\mathcal{P}({\mathbb{T}^{d}})$ are such
that $\mathcal{W}(\mu_{0},\mu_{1})<\infty$, then there exists a sequence
$\{(\mu_{t}^{n},\bm{\nu}_{t}^{n})_{t\in\lbrack0,1]}\}_{n\in\mathbb{N}}$ in
$\mathcal{CE}_{1}(\mu_{0},\mu_{1})$ such that
\[
\mathcal{W}(\mu_{0},\mu_{1})^{2}=\lim_{n\rightarrow\infty}\int_{0}%
^{1}\mathcal{A}(\mu_{t}^{n},\bm{\nu}_{t}^{n})\;dt.
\]
Thus, by Proposition \ref{compactness}, there exists a pair $(\mu
_{t},\bm{\nu}_{t})_{t\in\lbrack0,1]}\in\mathcal{CE}_{1}(\mu_{0},\mu_{1})$
such that, for some subsequence $(n_{k})$, we have $\mu_{t}^{n_{k}%
}\rightharpoonup\mu_{t}$ weakly in $\mathcal{P}({\mathbb{T}^{d}})$ for all
$t\in\lbrack0,T]$ and $\bm{\nu}^{n_{k}}\rightharpoonup\bm{\nu}$
weakly-$\ast$ in $\mathcal{M}_{\text{loc}}(G({\mathbb{T}^{d}})\times
\lbrack0,T])$. Therefore, also by the same proposition, we have
\[
\int_{0}^{1}\mathcal{A}(\mu_{t},\bm{\nu}_{t})\;dt\leqslant\liminf
_{k\rightarrow\infty}\int_{0}^{1}\mathcal{A}(\mu_{t}^{n_{k}},\bm{\nu}%
_{t}^{n_{k}})\;dt=\mathcal{W}(\mu_{0},\mu_{1})^{2}\;.
\]
Therefore
\[
\mathcal{W}(\mu_{0},\mu_{1})^{2}=\int_{0}^{1}\mathcal{A}(\mu_{t},\bm{\nu
}_{t})\;dt.
\]
Using Lemma \ref{WT2} and the Cauchy-Schwarz inequality, we conclude that
\[
\mathcal{W}(\mu_{0},\mu_{1})\leqslant\int_{0}^{1}\sqrt{\mathcal{A}(\mu
	_{t},\bm{\nu}_{t})}\;dt\leqslant\left(  \int_{0}^{1}\mathcal{A}(\mu
_{t},\bm{\nu}_{t})\;dt\right)  ^{\frac{1}{2}}=\mathcal{W}(\mu_{0},\mu
_{1}).
\]
This implies that $\mathcal{A}(\mu_{t},\bm{\nu}_{t})=\mathcal{W}(\mu
_{0},\mu_{1})$ for a.e. $t\in\lbrack0,1]$.

Moreover, if a pair $(\mu_{t},\bm{\nu}_{t})_{t\in\lbrack0,1]}%
\subset\mathcal{CE}_{1}(\mu_{0},\mu_{1})$ is such that $\mathcal{A}(\mu
_{t},\bm{\nu}_{t})$ is constant almost everywhere, we have that, for any
$0\leqslant s\leqslant r\leqslant1$, the pair $(\mu_{t+s},\bm{\nu}%
_{t+s})_{t\in\lbrack0,r-s]}$ belongs to $\mathcal{CE}_{r-s}(\mu_{s},\mu_{r})$
and, by Lemma \ref{WT2}, it satisfies
\[
\mathcal{W}(\mu_{s},\mu_{r})\leqslant\int_{0}^{r-s}\mathcal{A}(\mu
_{t+s},\bm{\nu}_{t+s})\;dt=(r-s)W(\mu_{0},\mu_{1})
\]
The opposite inequality is always true. In fact, if it is not verified for
some $s$ and $r$, then
\begin{align*}
	\mathcal{W}(\mu_{0},\mu_{1})  &  \leqslant\mathcal{W}(\mu_{0},\mu
	_{s})+\mathcal{W}(\mu_{s},\mu_{r})+\mathcal{W}(\mu_{r},\mu_{1})\\
	&  <s\mathcal{W}(\mu_{0},\mu_{1})+(r-s)\mathcal{W}(\mu_{0},\mu_{1}
	)+(1-r)\mathcal{W}(\mu_{0},\mu_{1})\\
	&  =\mathcal{W}(\mu_{0},\mu_{1}),
\end{align*}
which gives a contradiction. Therefore, one must have $\mathcal{W}(\mu_{s}%
,\mu_{r})=(r-s)\mathcal{W}(\mu_{0},\mu_{1})$. \fin

In next proposition we state the convexity of $\mathcal{W}^{2}$ in
$\mathcal{P}({\mathbb{T}^{d}}).$ The proof is similar to that of
\cite[Proposition 4.7]{erbar} and is left to the reader.

\begin{proposition}
	The function $(\mu_{0},\mu_{1})\in\mathcal{P}({\mathbb{T}^{d}})\times
	\mathcal{P}({\mathbb{T}^{d}})\mapsto\mathcal{W}(\mu_{0},\mu_{1})^{2}$ is convex.
\end{proposition}

\begin{definition}
	For $\mu_{0},\mu_{1}\in\mathcal{P}({\mathbb{T}^{d}})$ such that $\mathcal{W}%
	(\mu_{0},\mu_{1})<\infty$, we call any pair $(\mu_{t},\bm{\nu}_{t}%
	)_{t\in\lbrack0,1]}$ given by Proposition \ref{geodesic} as \emph{geodesic
		pair} or \emph{optimal pair} for $\mathcal{W}(\mu_{0},\mu_{1})$. And the
	curves $(\mu_{t})_{t\in\lbrack0,1]}$ as a geodesic curve between $\mu_{0}$ and
	$\mu_{1}$.
\end{definition}

\begin{proposition}
	The map $(\mu_{0},\mu_{1})\mapsto\mathcal{W}(\mu_{0},\mu_{1})$ is lower
	semicontinuous w.r.t. weak convergence. Moreover, the topology induced by
	$\mathcal{W}$ is stronger than the weak topology.
\end{proposition}

\emph{Proof.} Let $\{\mu_{0}^{n}\}_{n\in\mathbb{N}},\{\mu_{1}^{n}%
\}_{n\in\mathbb{N}}\subseteq\mathcal{P}({\mathbb{T}^{d}})$ be two sequences
such that $\mu_{i}^{n}\rightharpoonup\mu_{i}$ weakly for $i=0,1$ and $\mu
_{0},\mu_{1}\in\mathcal{P}({\mathbb{T}^{d}})$. Let us suppose that
$\liminf_{n\rightarrow\infty}\mathcal{W}(\mu_{0}^{n},\mu_{1}^{n})<\infty$ and
let $(n_{k})$ be a subsequence such that $\lim_{k\rightarrow\infty}%
\mathcal{W}(\mu_{0}^{n_{k}},\mu_{1}^{n_{k}})=\liminf_{n\rightarrow\infty
}\mathcal{W}(\mu_{0}^{n},\mu_{1}^{n})$. Proposition \ref{geodesic} assures
that for every $k$ there exists a pair $(\mu_{t}^{n_{k}},\bm{\nu}%
_{t}^{n_{k}})_{t\in\lbrack0,1]}\in\mathcal{CE}_{1}(\mu_{0}^{n_{k}},\mu
_{1}^{n_{k}})$ such that
\[
\mathcal{W}(\mu_{0}^{n_{k}},\mu_{1}^{n_{k}})^{2}=\int_{0}^{1}\mathcal{A}%
(\mu_{t}^{n_{k}},\bm{\nu}_{t}^{n_{k}})\;dt.
\]
Thus, the existence of the limit $\lim_{k}\mathcal{W}(\mu_{0}^{n_{k}},\mu
_{1}^{n_{k}})$ assures that the hypothesis of Proposition \ref{compactness} is
valid and we can extract a new subsequence, still denoted by $(\mu_{t}^{n_{k}%
},\bm{\nu}_{t}^{n_{k}})_{t\in\lbrack0,1]}$, such that $\mu_{t}^{n_{k}%
}\rightharpoonup\mu_{t}$ weakly in $\mathcal{P}({\mathbb{T}^{d}})$ for all
$t\in\lbrack0,T]$ and $\bm{\nu}^{n_{k}}\rightharpoonup\bm{\nu}$
weakly-$\ast$ in $\mathcal{M}_{\text{loc}}(G({[0,1)^{d}})\times\lbrack0,T])$
for some pair $(\mu_{t},\bm{\nu}_{t})_{t\in\lbrack0,1]}\in\mathcal{CE}%
_{1}(\mu_{0},\mu_{1}).$ Thus, from the lower semicontinuity result in the same
proposition, it follows that
\begin{align*}
	\mathcal{W}(\mu_{0},\mu_{1})^{2}  &  \leqslant\int_{0}^{1}\mathcal{A}(\mu
	_{t},\bm{\nu}_{t})\;dt\leqslant\liminf_{k\rightarrow\infty}\int_{0}
	^{1}\mathcal{A}(\mu_{t}^{n_{k}},\bm{\nu}_{t}^{n_{k}})\;dt=\lim
	_{k\rightarrow\infty}\mathcal{W}(\mu_{0}^{n_{k}},\mu_{1}^{n_{k}})^{2}\\
	&  =\liminf_{n\rightarrow\infty}\mathcal{W}(\mu_{0}^{n},\mu_{1}^{n})^{2}.
\end{align*}

Now, let $\{\mu_{n}\}_{n}$ be a sequence in $\mathcal{P}({\mathbb{T}^{d}})$
such that $\mathcal{W}(\mu_{n},\mu)\rightarrow0$. Then, any convergent
subsequence $\{\mu_{n_{k}}\}$ that converges weakly to a $\mu_{*}$ satisfies,
by the lower semicontinuity of $\mathcal{W}$, that $\mathcal{W}(\mu, \mu_{*})
\leqslant\liminf\mathcal{W}(\mu_{n_{k}},\mu)=0$. Since $\mathcal{P}%
({\mathbb{T}^{d}})$ is sequentially compact, then we obtain that $\mu_{n}$
converges weakly to $\mu$.

\fin

\begin{proposition}
	Given $\mu_{*}\in\mathcal{P}({\mathbb{T}^{d}})$, we have that the pair
	$\left(  \mathcal{P}_{\mu_{*}},\mathcal{W} \right)  $ is a complete metric
	space, where $\mathcal{P}_{\mu_{*}}:=\{\mu\in\mathcal{P}({\mathbb{T}^{d}})
	\;|\; \mathcal{W}(\mu_{*},\mu)<\infty\}$.
\end{proposition}

\emph{Proof.} Let $(\mu_{n})_{n\in\mathbb{N}}\subseteq\mathcal{P}_{\mu_{\ast}%
}$ be a Cauchy sequence w.r.t. $\mathcal{W}$. Since $\mathcal{P}%
({\mathbb{T}^{d}})$ is weakly compact, we can extract a subsequence
$(\mu_{n_{k}})_{k\in\mathbb{N}}$ such that $\mu_{n_{k}}\rightharpoonup\mu$ for
some $\mu\in\mathcal{P}({\mathbb{T}^{d}})$. Using the lower semicontinuity of
$\mathcal{W}$, we have
\begin{equation}
	\mathcal{W}(\mu,\mu_{m})\leqslant\liminf_{k\rightarrow\infty}\mathcal{W}
	(\mu_{n_{k}},\mu_{m})\;,\;\;\forall\;m\in\mathbb{N}. \label{cauchy}%
\end{equation}
Thus, since the right hand side of \eqref{cauchy} can be made arbitrarily
small if $m$ is sufficiently big, we conclude that $(\mu_{n})_{n\in\mathbb{N}%
}$ converges to $\mu,$ which therefore must belong to $\mathcal{P}_{\mu_{\ast
}}$. \fin


\section{Subdifferential Calculus}

\quad In this section we construct a subdifferential calculus in the space
$\mathcal{P}({\mathbb{T}^{d}})$ endowed with the non-local metric
$\mathcal{W}$. We give the corresponding notion of tangent space,
subdifferential of a functional and the characterization of this concept to
$\lambda-$convex functionals.

\begin{definition}
	For $\mu\in\mathcal{P}({\mathbb{T}^{d}})$ we define the tangent space of
	$\mathcal{P}({\mathbb{T}^{d}})$ at $\mu$ by
	\[
	\text{Tan}_{\mu}\mathcal{P}({\mathbb{T}^{d}}) := \left\{  \bm{\nu}
	\in\mathcal{M}_{\text{loc}}(G({\mathbb{T}^{d}})) \;\left|  \;
	\begin{array}
	[c]{l}%
	\bm{\cdot} \;\mathcal{A}(\mu,\bm{\nu})<\infty\;,\\
	\bm{\cdot} \; \mathcal{A}(\mu,\bm{\nu})\leqslant\mathcal{A}
	(\mu,\bm{\nu}+\eta)\;,\; \forall\eta\in\mathcal{M}_{\text{loc}
	}(G({\mathbb{T}^{d}})) \mbox{ with } \overline{\nabla}\cdot\eta= 0
	\end{array}
	\right.  \right\}
	\]
	
\end{definition}

If $d\mu(x) = \rho(x) dx$ then%
\[
\text{Tan}_{\mu}\mathcal{P}({\mathbb{T}^{d}}) := \left\{  w\widehat{\rho}%
_{m}\;dK^{\sigma}\in\mathcal{M}_{\text{loc}}(G({\mathbb{T}^{d}})) \;|\;
\exists(\phi_{n})_{n\in\mathbb{N}}\subset C_{c}^{\infty}({\mathbb{T}^{d}%
})\;,\;\overline{\nabla}\phi_{n} \rightarrow w \;\; in \;L^{2}(\widehat{\rho
}_{m}dK^{\sigma}) \right\}
\]

Let $\mu_{0},\mu_{1}\in\mathcal{P}({\mathbb{T}^{d}})$ be such that
$\mathcal{W}(\mu_{0},\mu_{1})<\infty$, then by Proposition \ref{geodesic} there exists a geodesic pair
$(\mu_{t},\bm{\nu}_{t})_{t\in\lbrack0,1]}\in\mathcal{CE}_{1}(\mu_{0}%
,\mu_{1})$ satisfying
\[
\mathcal{W}^{2}(\mu_{0},\mu_{1})=\mathcal{A}(\mu_{t},\bm{\nu}%
_{t})\;,\;\;\mbox{
	for a.e. }t\in\lbrack0,1].
\]
If $K\subseteq G({\mathbb{T}^{d}})$ is compact, then by Lemma
\ref{comp.estimate} there exists a constant $C=C(K)$ such that the total
variation of $\bm{\nu}_{t}$ is uniformly bounded for a.e. $t\in
\lbrack0,1]$ with
\[
|\bm{\nu}_{t}|(K)\leqslant C(K)\sqrt{\mathcal{A}(\mu_{t},\bm{\nu}
	_{t})}=C(K)\mathcal{W}(\mu_{0},\mu_{1}).
\]
Therefore, there exists a sequence $t_{n}\rightarrow0^{+}$ and $\bm{\nu
}\in\mathcal{M}_{\text{loc}}(G({[0,1)^{d}}))$ such that $\bm{\nu}_{t_{n}%
}\overset{\ast}{\rightharpoonup}\bm{\nu}$. Thus, since the set of weak-$* $ limit measures for $\bm{\nu}_{t}$ as $t\to 0$ is non empty, we can have the following definition:

\begin{definition}
	Given $\mu_{0},\mu_{1}\in\mathcal{P}({\mathbb{T}^{d}})$ such that
	$\mathcal{W}(\mu_{0},\mu_{1})< \infty$, let $\Gamma_{0}(\mu_{0},\mu_{1})$ be
	the set of all weak-$* $ limit measures of $\bm{\nu}_{t}$ as $t\rightarrow
	0^{+}$ where $(\mu_{t},\bm{\nu}_{t})_{t \in[0,1] } \in\mathcal{CE}_{1}
	(\mu_{0},\mu_{1})$ is a geodesic pair.
\end{definition}

Now we give the notion of subdifferential for functionals defined in
$\mathcal{P}({\mathbb{T}^{d}})$ with respect to the geometric structure given
by $\mathcal{W}$.

\begin{definition}
	Let $\mu\in\mathcal{P}({\mathbb{T}^{d}})$ with $d\mu(x)=\rho(x)dx$ and let
	$\mathcal{E}:\mathcal{P}({\mathbb{T}^{d}})\rightarrow(-\infty,+\infty]$ be a
	proper and weak lower semicontinuous functional. We say that $d\xi
	=\zeta\widehat{\rho}_{m}dK^{\sigma}\in\text{Tan}_{\mu}\mathcal{P}%
	({\mathbb{T}^{d}})$ belongs to the subdifferential $\partial\mathcal{E}(\mu)$
	if
	\[
	\liminf_{\mu_{1}\rightarrow\mu}\frac{1}{\mathcal{W}(\mu,\mu_{1})}\left(
	\mathcal{E}(\mu_{1})-\mathcal{E}(\mu)-\inf_{\bm{\nu}\in\Gamma_{0}(\mu
		,\mu_{1})}\int_{G({[0,1)^{d}})}\zeta(x,y)\;d\bm{\nu}(x,y)\right)
	\geqslant0
	\]
	or, equivalently,
	\[
	\mathcal{E}(\mu_{1})-\mathcal{E}(\mu)\geqslant\inf_{\bm{\nu}\in\Gamma
		_{0}(\mu,\mu_{1})}\int_{G({[0,1)^{d}})}\zeta(x,y)\;d\bm{\nu}%
	(x,y)+o(\mathcal{W}(\mu,\mu_{1})).
	\]
	
\end{definition}

The next lemma deals with geodesics constructed via a glueing type process.
See \cite{DNS} for such a type of construction in another context.

\begin{lemma}
	\label{glueing} Let $(\mu_{t},\bm{\nu}_{t})_{t\in\lbrack0,1]}%
	\in\mathcal{CE}_{1}(\mu_{0},\mu_{1})$ be a geodesic pair for $\mathcal{W}%
	(\mu_{0},\mu_{1})$ and let $s\in(0,1)$. If $(\mu_{t}^{s},\bm{\nu}_{t}%
	^{s})_{t\in\lbrack0,s]}$ is another geodesic pair for $\mathcal{W}(\mu_{0}%
	,\mu_{s})$, then the composition
	\[
	\overline{\mu}_{t}:=\left\{
	\begin{array}
	[c]{cc}%
	\mu_{t}^{s}\;, & 0\leqslant t\leqslant s\\
	\mu_{t}\;, & s<t\leqslant1\\
	&
	\end{array}
	\right.  \overline{\bm{\nu}}_{t}:=\left\{
	\begin{array}
	[c]{cc}%
	\bm{\nu}_{t}^{s}\;, & 0\leqslant t\leqslant s\\
	\bm{\nu}_{t}\;, & s<t\leqslant1\\
	&
	\end{array}
	\right.
	\]
	is also a geodesic pair for $\mathcal{W}(\mu_{0},\mu_{1})$.
\end{lemma}

Let us define what a $\lambda$-geodesically convex functional means in our context. We say that a proper functional $\mathcal{E}:\mathcal{P}({\mathbb{T}^{d}})\rightarrow(-\infty,+\infty]$ is $\lambda$-convex for some $\lambda\in \R$ if for every $\mu_0,\mu_1\in \PT$ with $\W(\mu_0,\mu_1)<\infty$, there exists a geodesic curve $(\mu_t)_{t\in[0,1]}$ such that the following inequality is true
\[\mathcal{E}(\mu_t)\meni (1-t)\mathcal{E}(\mu_0)+t\mathcal{E}(\mu_1)-\frac{\lambda}{2}t(1-t)\W(\mu_0,\mu_1)^2,\;\;\forall \;t \in [0,1].\]

In the case when the functional $\mathcal{E}$ is $\lambda$-geodesically convex, we have the
following characterization.

\begin{proposition}
	Let $\mathcal{E}:\mathcal{P}({\mathbb{T}^{d}})\rightarrow(-\infty,+\infty]$ be
	a lower semicontinuous and $\lambda$-geodesically convex functional. If $d\mu(x)=\rho(x)dx$
	then $d\xi=\zeta\widehat{\rho}_{m}dK^{\sigma}\in\partial\mathcal{E}(\mu),$ if
	and only if, for all $\mu_{1}\in\mathcal{P}({\mathbb{T}^{d}})$ such that
	$\mathcal{W}(\mu,\mu_{1})<\infty$ there exists $\bm{\nu}_{0}\in\Gamma
	_{0}(\mu,\mu_{1})$ satisfying
	\[
	\mathcal{E}(\mu_{1})-\mathcal{E}(\mu)\geqslant\int_{G({[0,1)^{d}})}%
	\zeta(x,y)\;d\bm{\nu}_{0}(x,y)+\frac{\lambda}{2}\mathcal{W}^{2}(\mu
	,\mu_{1}).
	\]
	
\end{proposition}

\emph{Proof.} Let $\mu_{1}\in\mathcal{P}({\mathbb{T}^{d}})$ be such that
$\mathcal{W}(\mu,\mu_{1})<\infty$ and let $(\mu_{t},\bm{\nu}_{t}%
)_{t\in\lbrack0,1]}\in\mathcal{CE}_{1}(\mu,\mu_{1})$ be a geodesic pair for
$\mathcal{W}(\mu,\mu_{1})$. Since $\mathcal{W}(\mu,\mu_{t})=t\mathcal{W}%
(\mu,\mu_{1})$ for all $t\in\lbrack0,1]$, there exists $\bm{\nu}_{0}%
^{t}\in\Gamma_{0}(\mu,\mu_{t})$ such that
\begin{equation}
	\mathcal{E}(\mu_{t})-\mathcal{E}(\mu)\geqslant\int_{G({[0,1)^{d}})}
	\zeta(x,y)\;d\bm{\nu}_{0}^{t}+o(t). \label{o(t)}%
\end{equation}
\textit{Claim:} $t^{-1}\bm{\nu}_{0}^{t}\in\Gamma_{0}(\mu,\mu_{1}).$

Let $(\mu_{s}^{t},\bm{\nu}_{s}^{t})_{s\in\lbrack0,1]}\in\mathcal{CE}%
_{1}(\mu,\mu_{t})$ be a geodesic pair for $\mathcal{W}(\mu,\mu_{t})$ such that
$\bm{\nu}_{0}^{t}$ is a limit point of $(\bm{\nu}_{s}^{t})_{s}$ as
$s\rightarrow0^{+}$. By Lemma \ref{rescaling}, the pair $\left(  \mu_{\frac
	{s}{t}}^{t},t^{-1}\bm{\nu}_{\frac{s}{t}}^{t}\right)  _{s\in\lbrack0,t]}$
are still geodesic pairs for $\mathcal{W}(\mu,\mu_{t})$ for all $t$ and, by
Lemma \ref{glueing}, the pairs $(\overline{\mu}^{t}_{s},\overline{\bm{\nu
	}}^{t}_{s})_{s\in\lbrack0,1]}$ defined by
	\[
	\overline{\mu}_{s}^{t}:=\left\{
	\begin{array}
	[c]{cc}%
	\mu_{\frac{s}{t}}^{t}\;, & 0\leqslant s\leqslant t\\
	\mu_{s}\;, & t<s\leqslant1
	\end{array}
	\right.  \overline{\bm{\nu}}_{s}^{t}:=\left\{
	\begin{array}
	[c]{cc}%
	t^{-1}\bm{\nu}_{\frac{s}{t}}^{t}\;, & 0\leqslant s\leqslant t\\
	\bm{\nu}_{s}\;, & t<s\leqslant1
	\end{array}
	\right.
	\]
	are also geodesic pairs for $\mathcal{W}(\mu,\mu_{1})$ and have $t^{-1}%
	\bm{\nu}_{0}^{t}$ as a limit point for $\overline{\bm{\nu}}_{s}$ as
	$s\rightarrow0^{+}.$
	
	Now, by Proposition \ref{geodesic}, these geodesic pairs satisfy $W^{2}%
	(\mu,\mu_{1})=\mathcal{A}(\overline{\mu}_{s}^{t},\overline{\bm{\nu}}%
	_{s}^{t})$ for all $t\in\lbrack0,1]$ and for a.e. $s\in\lbrack0,1]$, and for
	all $K\subseteq G({\mathbb{T}^{d}})$ we can use Lemma \ref{comp.estimate} and
	the weak-weak$\ast$ lower semicontinuity of $\mathcal{A}$ to estimate
	\[
	t^{-1}|\bm{\nu}_{0}^{t}|(K)\leqslant C(K)\sqrt{\mathcal{A}(\mu
		,t^{-1}\bm{\nu}_{0}^{t})}\leqslant C(K)\liminf_{s\rightarrow0}%
	\sqrt{\mathcal{A}(\overline{\mu}_{s}^{t},\overline{\bm{\nu}}_{s}^{t}
		)}=C(K)\mathcal{W}(\mu,\mu_{1}).
	\]
	Therefore, there exists $\bm{\nu}_{0}$ weak-$\ast$ limit point of
	$t^{-1}\bm{\nu}_{0}^{t}$ as $t\rightarrow0^{+}$. Finally, by Proposition
	\ref{compactness}, there exists also a geodesic pair $(\overline{\overline
		{\mu}}_{t},\overline{\overline{\bm{\nu}}}_{t})_{t\in\lbrack0,1]}%
	\in\mathcal{CE}_{1}(\mu,\mu_{1})$ such that $\overline{\overline{\bm{\nu}%
		}}_{t}\overset{\ast}{\rightharpoonup}\bm{\nu}_{0}$ and therefore,
		$\bm{\nu}_{0}\in\Gamma_{0}(\mu,\mu_{0})$.
		
		Since $\mathcal{E}$ is $\lambda$-convex, the following inequality holds for
		all $t\in(0,1]$
		\begin{equation}
			\label{eqlambda1}\mathcal{E}(\mu_{1})-\mathcal{E}(\mu)-\frac{\lambda}%
			{2}(1-t)\mathcal{W} ^{2}(\mu,\mu_{1})\geqslant\frac{\mathcal{E}(\mu
				_{t})-\mathcal{E}(\mu)}{t}%
		\end{equation}
		and, for $d\xi=\zeta\widehat{\rho}_{m}dK^{\sigma}\in\partial\mathcal{E}(\mu)$,
		we obtain from \eqref{o(t)} that
		\begin{equation}
			\label{eqlambda2}\liminf_{t\rightarrow0}\frac{\mathcal{E}(\mu_{t}%
				)-\mathcal{E}(\mu)} {t}\geqslant\liminf_{t\rightarrow0}\int_{G({[0,1)^{d}}%
				)}\zeta(x,y)t^{-1} \;d\bm{\nu}_{0}^{t}=\int_{G({[0,1)^{d}})}%
			\zeta(x,y)\;d\bm{\nu} _{0}\text{,}%
		\end{equation}
		where the equality follows from the following fact: since $d\mu(x)=\rho
		(x)\;dx$ and $\mathcal{A}(\mu,t^{-1}\bm{\nu}_{0}^{t}),\mathcal{A}%
		(\mu,\bm{\nu}_{0})<\infty$ then there exist $w_{0}^{t},w_{0}\in
		L^{2}(\widehat{\rho}_{m}\;dK^{\sigma})$ such that $t^{-1}\;d\bm{\nu}%
		_{0}^{t}=w_{0}^{t}\widehat{\rho}_{m}\;dK^{\sigma}$ and $d\bm{\nu}%
		_{0}=w_{0}\widehat{\rho}_{m}\;dK^{\sigma}$. Therefore, $t^{-1}\bm{\nu}%
		_{0}^{t}\overset{\ast}{\rightharpoonup}\bm{\nu}_{0}$ implies that
		$w_{0}^{t}\rightharpoonup w_{0}$ in $L^{2}(\widehat{\rho}_{m}\;dK^{\sigma})$.
		Thus, since $\zeta\in L^{2}(\widehat{\rho}_{m}\;dK^{\sigma})$, we have
		\begin{align*}
			\lim_{t\rightarrow0}\int_{G({[0,1)^{d}})}\zeta(x,y)\;t^{-1}\;d\bm{\nu}
			_{0}^{t}  &  =\lim_{t\rightarrow0}\int_{G({[0,1)^{d}})}\zeta(x,y)w_{0}
			^{t}(x,y)\widehat{\rho}_{m}\;dK^{\sigma}\\
			&  =\int_{G({[0,1)^{d}})}\zeta(x,y)w_{0}(x,y)\widehat{\rho}_{m}\;dK^{\sigma}\\
			&  =\int_{G({[0,1)^{d}})}\zeta(x,y)\;d\bm{\nu}_{0}.
		\end{align*}
		Finally, from \eqref{eqlambda1} and \eqref{eqlambda2} we obtain the desired expression.
		
		\fin
		
		\begin{remark}
			If $\Omega\subset\mathbb{R}^{k}$ is open and $\bm{\nu}_{n}\overset{\ast
			}{\rightharpoonup}\bm{\nu}$ in $\mathcal{M}_{\text{loc}}(\Omega)$ with
			$d\bm{\nu}_{n}=w_{n}d\bm{\pi}$ for some $0\leqslant\pi\in
			\mathcal{M}_{\text{loc}}(\Omega)$, then $w_{n}\rightharpoonup w$ in
			$L^{2}(d\bm{\pi})$.
		\end{remark}
	

\section{Generalized minimizing movements}
		
\quad In this section we consider the R\'{e}nyi entropy functional which is the
internal energy functional associated to the porous medium equation in the
classical case. We show that for each $\mu_{0}\in\mathcal{P}({\mathbb{T}^{d}})$
there exists an absolutely continuous curve obtained as the limit (up to a
subsequence) of the discrete variational scheme \eqref{eq:discr-schem}. As
commented in the introduction, such curve is a generalized minimizing movement
for \eqref{FPMeq} starting at $\mu_{0}$.
		
\begin{definition}
			\label{def:intern-energ} For $\mu\in\mathcal{P}({\mathbb{T}^{d}})$ and
			$m\in(0,2]$ we define the entropy $\mathcal{U}_{m}$ at $\mu$ by
			\[
			\mathcal{U}_{m}(\mu):=\left\{
			\begin{array}
			[c]{ll}%
			\displaystyle\int_{{[0,1)^{d}}}U_{m}(\rho(x))\;dx\;, & \mbox{ if }d\mu
			=\rho(x)\;dx\\
			+\infty\;, & \mbox{ else }\\
			&
			\end{array}
			\right.  ,
			\]
			where
			\[
			U_{m}(s):=\left\{
			\begin{array}
			[c]{ll}%
			s\log s & \mbox{ if }\;m=1\\
			\displaystyle\frac{s^{m}}{m-1} & \mbox{ else}\\
			&
			\end{array}
			\right.  .
			\]
			
		\end{definition}
		
		\begin{definition}
			For $\mu\in\mathcal{P}({\mathbb{T}^{d}})$ and $m=1$ we define the entropy
			dissipation (or fisher information associated) $\mathcal{I}_{1}$ at $\mu$ by
			\[
			\mathcal{I}_{1}(\mu):= \Xi(\bm{\mu}^{1},\bm{\mu}^{2}):=\frac{1}{2}
			\int_{G({[0,1)^{d}})}\frac{(\rho^{1}-\rho^{2})^{2}}{\theta(\rho^{1},\rho^{2}
				)}\;d\bm{\lambda}
			\]
			where $d\bm{\mu}^{1}(x,y)=K^{\sigma}(x-y)dyd\mu(x)$, $d\bm{\mu}
			^{2}(x,y)=K^{\sigma}(x-y)dxd\mu(y)$ and $\lambda\in\mathcal{M}_{\text{loc}
			}(G({\mathbb{T}^{d}}))$ is such that $\bm{\mu}^{1},\bm{\mu}^{2}
			\ll\lambda$.
		\end{definition}
		
		\begin{remark}
			Note that when $d\mu(x)=\rho(x)dx$, we have
			\[
			\mathcal{I}_{1}(\mu)=\frac{1}{2}\int_{G({[0,1)^{d}})}\overline{\nabla}
			\rho(x,y)\overline{\nabla}\log{\rho}(x,y)K^{\sigma}(x-y)dxdy.
			\]
			Since the integrand in the definition of $\Xi$ is a lower semicontinuous
			function of $\rho^{1}$ and $\rho^{2}$, and also positively 1-homogeneous, we
			have that $\Xi$ is lower semicontinuous w.r.t. the weak* convergence in
			$\mathcal{M}_{\text{loc}}(G({\mathbb{T}^{d}}))$. Therefore, if $\mu
			_{n}\rightharpoonup\mu$ weakly in $\mathcal{P}({\mathbb{T}^{d}})$, then
			$\bm{\mu}_{n}^{i}\overset{\ast}{\rightharpoonup}\bm{\mu}^{i}$, for
			$i=1,2$, and it implies that
			\[
			\mathcal{I}_{1}(\mu)=\Xi(\bm{\mu}^{1},\bm{\mu}^{2})\leqslant
			\liminf_{n\rightarrow\infty}\Xi(\bm{\mu}_{n}^{1},\bm{\mu}_{n}
			^{2})=\liminf_{n\rightarrow\infty}\mathcal{I}_{1}(\mu_{n}).
			\]
			Analogously, for $m\neq1$ we can define the fisher information as
			\[
			\mathcal{I}_{m}(\mu):=\frac{m}{m-1}\int_{G([0,1)^{d})}{\overline{\nabla}
				(\rho^{m-1})\overline{\nabla}(\rho^{m})K^{\sigma}(x-y)\ dxdy},
			\]
			for $d\mu=\rho\ dx$. Note that, in this case the integrand is $(2m-1)-$
			homogeneous and thus it is not possible to extend it for more general measures.
		\end{remark}
		
		The next theorem shows a coercivity property to the internal energy defined in
		Definition \ref{def:intern-energ}.
		
		\begin{theorem}
			\label{Phi.bounded} For any $\tau>0$ and $\mu_{*} \in\mathcal{P}
			({\mathbb{T}^{d}})$, the function
			\[
			\mu\in\mathcal{P}_{\mu_{*}} \mapsto\Phi(\tau,\mu_{*};\mu) := \frac{1}{2\tau
			}\mathcal{W}^{2}(\mu_{*},\mu) + \mathcal{U}_{m}(\mu)
			\]
			is bounded from below in $\mathcal{P}({\mathbb{T}^{d}})$. Moreover, there
			exists a unique $\mu^{0} \in\mathcal{P}({\mathbb{T}^{d}})$ (depending on
			$\tau$ and $\mu_{*}$) such that
			\[
			\Phi(\tau,\mu_{*};\mu^{0}) \leqslant\Phi(\tau,\mu_{*};\mu)\;,\;\; \forall
			\;\mu\in\mathcal{P}({\mathbb{T}^{d}}).
			\]
			
		\end{theorem}
		
		\emph{Proof.} We first notice that for $m=1$ and any $\mu= \rho(x)\;dx
		\in\mathcal{P}_{ac}({\mathbb{T}^{d}})$, we can rewrite $\mathcal{U}_{1}$ as
		\[
		\mathcal{U}_{1}(\mu) = \int_{[0,1)^{d}}\rho\log\rho\;dx = \int_{[0,1)^{d}%
		}\left(  \rho\log\rho-\rho+1\right)  \;dx
		\]
		and use the fact that $s\mapsto s \log s -s +1$ is nonnegative in $[0,\infty).$
		
		If $m\in(1,2]$ then the nonnegativity of $\mathcal{U}_{m}$ is clear.
		
		Finally, if $m \in(0,1)$, we have that
		\[
		\mathcal{U}_{m}(\mu) = \int_{{[0,1)^{d}}}\left(  \frac{\rho(x)^{m}}{m-1}%
		+\rho(x)-1 \right)  \;dx
		\]
		and we can use the fact that $s\mapsto\frac{s^{m}}{m-1}+s-1$ is bounded from below.
		
		Therefore, for all $m\in(0,2]$, $\tau>0$ and $\mu_{*}\in\mathcal{P}%
		({\mathbb{T}^{d}})$, the functional $\Phi$ is bounded from below.
		
		The existence of points of minimum follows from the fact that any minimizing
		sequence for $\Phi$ is relatively compact w.r.t. the weak convergence.
		Therefore, by the lower semicontinuity of $\mathcal{W}$ and $\mathcal{U}_{m}$ (the semicontinuity of the functional w.r.t. the weak convergence is a very well known fact),
		we have that any limit point attains the minimum value for $\Phi$. The
		uniqueness is a consequence of the fact that the composition $t\in
		[0,1]\mapsto\mu_{t}:=(1-t)\mu_{0}+t\mu_{1}\mapsto\mathcal{U}_{m}(\mu_{t})$ is
		strictly convex. \fin
		
		\begin{theorem}
			\label{GMMcurve} For any $\tau>0$ and $\mu_{0}\in\mathcal{P}({\mathbb{T}^{d}%
			})$ such that $\mathcal{U}_{m}(\mu_{0})<\infty$ we can define
			\[
			\mu_{\tau}^{0}:=\mu_{0}\;\;\mbox{ and }\;\;\mu_{\tau}^{n}:=\emph{argmin}%
			\left\{  \left.  \Phi(\tau,\mu_{0}^{n-1};\mu)\;\right\vert \;\mu\in
			\mathcal{P}({\mathbb{T}^{d}})\right\}  \;,\;\;\forall n\in\mathbb{N}.
			\]
			Now defining the interpolation $\mu_{\tau}:[0,\infty)\rightarrow
			\mathcal{P}({\mathbb{T}^{d}})$ by
			\[
			\mu_{\tau}(t):=\mu_{\tau}^{n}\;\;\;\mbox{ for }\;\;t\in\lbrack n\tau
			,(n+1)\tau)\;,\;\;\mbox{ and }\;n\in\mathbb{N}\cup\{0\},
			\]
			we have that there exists a curve $\mu\in AC_{\emph{loc}}([0,\infty
			),\mathcal{P}({\mathbb{T}^{d}}))$ such that (up to a subsequence)
			\[
			\mu_{\tau}(t)\rightharpoonup\mu(t),\;\mbox{ as }\;\tau\rightarrow
			0\;\;\forall\;t\geqslant0.
			\]
			
		\end{theorem}
		
		\emph{Proof.} For a fixed $\mu_{*}\in\mathcal{P}({\mathbb{T}^{d}})$, we have
		from the previous results that the following ones are true:
		
		\begin{enumerate}
			\item The weak topology on $\mathcal{P}_{\mu_{*}}$ is compatible with
			$\mathcal{W}$, in the sense that it is weaker than the one generated by
			$\mathcal{W}$ and the function $(\mu_0,\mu_1)\mapsto \mathcal{W}(\mu_0,\mu_1)$ is weak lower semicontinuous;
			
			\item The entropy functional is weak lower semicontinuous on $\mathcal{P}
			_{\mu_{*}}$;
			
			\item For all $\tau>0$ the functional $\mu\in\mathcal{P}_{\mu_{\ast}}
			\mapsto\Phi(\tau,\mu_{\ast};\mu)$ is bounded from below;
			
			\item The space $\mathcal{P}_{\mu_{*}}$ is weakly sequentially compact.
		\end{enumerate}
		
		Therefore, we just need to apply the general results for metric spaces like in
		\cite[Proposition 2.2.3 or Corollary 3.3.4]{ambrosiogiglisavare}.
		
		\fin

%
%
%

		\noindent\textbf{Acknowledgement:} LCFF was supported by CNPq and FAPESP,
		Brazil. MCS acknowledges the support from the S\~{a}o Paulo Research
		Foundation (FAPESP) grant $\#$2014/23326-1, Brazil. JCV-G acknowledges the
		support from CNPq, Brazil.

\begin{small}

	\end{small}
	\end{document}